\documentclass[journal]{IEEEtran}

\usepackage{cite}
\usepackage{url}
\usepackage{graphicx}
\usepackage{color}
\usepackage{bm}
\usepackage{cite}
\usepackage{amsmath,amssymb,amsfonts}
\usepackage{algorithmic}
\usepackage{graphicx}
\usepackage{textcomp}
\usepackage{algorithmic}
\usepackage{algorithm}
\usepackage{array}
\usepackage[caption=false,font=normalsize,labelfont=sf,textfont=sf]{subfig}
\usepackage{stfloats}
\usepackage{url}
\usepackage{verbatim}
\usepackage{enumerate}
\usepackage{tikz}
\usepackage{bm}
\usepackage{threeparttable,booktabs}
\usepackage{blkarray}
\usepackage{amsmath}
\usepackage{multirow}

\usepackage{amsmath}
\usepackage{amsthm,amsfonts}
\usepackage{pgfplots}

\usepackage{tikz}
\usepackage[customcolors]{hf-tikz}

\usepgfplotslibrary{colorbrewer}
\usetikzlibrary{%
  3d,
  arrows,%
  shapes.misc,
  shapes.arrows,%
  chains,%
  matrix,%
	fit,%
  positioning,
  scopes,%
  decorations.pathmorphing,
  shadows,%
  decorations.markings,
  shapes.geometric,
  arrows.meta,bending,
	patterns,
	intersections, 
	pgfplots.fillbetween
}

\usepackage{tikz-3dplot}
\usepackage[pdfusetitle]{hyperref}
\hypersetup{
    colorlinks,
    linkcolor={red!50!black}, 
    citecolor={blue!50!black},
    urlcolor={blue!50!black}
}

\usepackage{calc}  
\newlength{\nomlabelwidth}                   
\newcommand{\setnomlabel}[1]{%
  \settowidth{\nomlabelwidth}{\textit{#1}}
}
\newenvironment{Nomenclature}[1]{%
  \setnomlabel{#1}%
  \section*{Nomenclature}
  \begin{tabular}{@{}p{\nomlabelwidth}@{\hspace{1em}}p{\linewidth-\nomlabelwidth-1em}@{}}
}{%
  \end{tabular}
}

\usepackage{cellspace}
\providecommand*{\mat}[1]{\mathbf#1}
\providecommand*{\M}[1]{\mathbf#1}
\providecommand*{\tM}[1]{\tilde{\mathbf#1}}

\providecommand*{\mrm}[1]{\mathrm{#1}}

\renewcommand{\vec}[1]{{\boldsymbol#1}}

\DeclareMathAccent{\ring}{\mathalpha}{operators}{"17}

\renewcommand{\j}{\mathrm{j}}

\newcommand{\ie}{\textit{i.e.}\/, }
\newcommand{\eg}{\textit{e.g.}\/, }

\newcommand{\cbullet}{\raisebox{0.25ex}{\textbullet}} 

\colorlet{dpurple}{blue!50!red}
\colorlet{dblue}{blue!50!black}
\colorlet{dgreen}{green!50!black}
\colorlet{dred}{red!50!black}
\colorlet{dyellow}{yellow!50!black}
\colorlet{dorange}{orange!50!black}
\definecolor{metal}{RGB}{218,165,32}
\definecolor{diel}{RGB}{1,165,32}
\definecolor{antenna}{RGB}{100,150,162}
\definecolor{breg}{rgb}{0.2,0.6,0.8}%
\definecolor{preg}{rgb}{0.8,0.2,0.2}%
\definecolor{reg}{RGB}{218,165,32}

\graphicspath{{figures/}{tikz/}}

\makeatother

\begin{document}
\title{Generalized Scattering Matrix Framework for Modeling Implantable Antennas in Multilayered Spherical Media}

\author{Chenbo Shi, Xin Gu, Shichen Liang, Jin Pan
\thanks{Manuscript received Jul. 17, 2025. Revised Nov. 21, 2025. (\textit{Corresponding author: Jin Pan.})}
\thanks{Chenbo Shi, Xin Gu, Shichen Liang and Jin Pan are with the School of Electronic Science and Engineering, University of Electronic Science and Technology of China, Chengdu 611731 China  (e-mail: chenbo\_shi@163.com; xin\_gu04@163.com; lscstu001@163.com; panjin@uestc.edu.cn).}
}


\maketitle

\begin{abstract}
This paper presents a unified and computationally efficient framework for modeling antennas embedded in spherically stratified media, applicable to implantable biomedical antennas and radome-enclosed systems. The method separates the characterization of the radiator from that of the surrounding medium by combining the antenna's free-space generalized scattering matrix (GSM) with a set of extended spherical scattering operators (SSOs). This decoupling enables rapid reevaluation under arbitrary changes of the spherical medium without re-simulating the antenna, yielding orders-of-magnitude speedups over traditional DGF-based MoM approaches. The SSO formulation accommodates multilayer, radially inhomogeneous, and radially uniaxial anisotropic profiles, and the GSM can be obtained from diverse numerical solvers or far-field data, supporting array-level synthesis and measurement-driven modeling. Extensive examples confirm excellent agreement with full-wave and DGF-based solutions, demonstrating the accuracy, generality, and practical versatility of the proposed framework.
\end{abstract}

\begin{IEEEkeywords}
  Generalized scattering matrix, implantable antennas, spherical layers, dyadic Green's functions.
\end{IEEEkeywords}

\begin{Nomenclature}{$\mat{\Gamma}$, $\mat{R}, \mat{T}, \mat{S}$}  
$\tau\sigma m l$ & Multiple indices of spherical modes.\\
$n$, $\tau n$ & Combined index, both representing $\tau\sigma m l$.\\
$\vec A_{\tau n}$ & Vector spherical harmonics.\\
$\mat{a}^\mrm{b},\mat{f}^\mrm{b}$ & Expansion vector of $\vec{u}_{n}^{\left( 1 \right)}, \vec{u}_{n}^{\left( 4 \right)}$ in bubble.\\
$a^\mrm{b}_{\tau n},f^\mrm{b}_{\tau n}$ & $\tau n$ entries of $\mat{a}^\mrm{b},\mat{f}^\mrm{b}$.\\
$\mat{a}^\mrm{f},\mat{f}^\mrm{f}$ & Expansion vector of $\vec{u}_{n}^{\left( 1 \right)}, \vec{u}_{n}^{\left( 4 \right)}$ in free-space.\\
$a^\mrm{f}_{\tau n},f^\mrm{f}_{\tau n}$ & $\tau n$ entries of $\mat{a}^\mrm{f},\mat{f}^\mrm{f}$.\\
$\zeta_l, \zeta_l^\prime, \xi_l, \xi_l^\prime$ & Riccati-Bessel functions and their derivatives.\\
$\mat{t}, \vec{\rho}$ & Spherical transition and reflection operators.\\
$t_{\tau n}, \rho_{\tau n}$ & $\tau n$ entries of $\mat{t}, \vec{\rho}$.\\
$\mat{\Psi}, \mat{\Phi}$ & Spherical bidirectional propagation operators.\\
$\Psi_{\tau n}, \Phi_{\tau n}$ & $\tau n$ entries of $\mat{\Psi}, \mat{\Phi}$.\\
$\mat{\Gamma}, \mat{R}, \mat{T}, \mat{S}$ & Blocks of generalized scattering matrix.\\
$R_{\tau n}^{\left(p\right)}(kr)$ & Auxiliary functions defined in Appendix \ref{APP_A}.\\
$x_\mrm{b}$, $x_\mrm{a}$ & $x_\mrm{b}=k_\mrm{b} r_\mrm{b},x_\mrm{a} = k_\mrm{f}r_\mrm{a}$.\\
\end{Nomenclature}

\section{Introduction}

\IEEEPARstart{I}{mplantable} antennas have gained broad interest owing to their importance in in-body wireless links, health monitoring, and telemedicine. In many works, the electromagnetic environment surrounding an implantable antenna—typically biological tissue—is modeled as a spherically stratified medium \cite{ref_implant1,ref_implant2,ref_implant3,ref_implant4,ref_implant5,ref_implant6}. Such models also arise in non-biological scenarios, including antennas enclosed by engineered radomes or metamaterial shells exhibiting inhomogeneous, lossy, or anisotropic properties. Throughout this paper, the term ``implantable'' antenna is therefore used in this generalized sense to denote radiators embedded within spherical-layered media.

Direct full-wave simulation of antennas in these environments is computationally demanding. A common alternative leverages closed-form dyadic Green’s functions (DGFs) for spherical stratifications in combination with the method of moments (MoM) \cite{ref_DGF1,ref_DGF2,ref_DGF3,ref_DGF4,ref_DGF5}, or employs hybrid T-matrix/MoM formulations that avoid explicit DGF construction \cite{ref_Capek1,ref_Capek2}. Although more efficient than brute-force discretization, these methods remain fundamentally tied to DGFs. Any modification of layer permittivity, permeability, or thickness requires recomputing the MoM impedance matrix, making parametric studies and optimization costly. The matrix dimension also scales with geometric complexity, limiting applicability to simple antennas. Furthermore, the utilization of DGFs are inseparably bound to the MoM framework, constraining compatibility with numerical schemes beyond MoM.

To address these challenges, this work introduces a separated modeling framework based on the generalized scattering matrix (GSM) of the antenna \cite{ref_scattering_Montgomery} and a set of spherical scattering operators (SSOs) representing the surrounding medium. The key idea is a clean decoupling: the antenna is fully characterized by its GSM, while all effects of the spherical stratification are captured by the SSOs. Because the GSM is constructed from a global spherical-wave expansion, it is independent of geometric detail and naturally provides model-order reduction. Once the GSM is obtained, changes to the enclosing layers require updating only the SSOs, not re-simulating the antenna. This separation yields significant computational savings and enables efficient design and optimization of spherical-layered structures such as radomes and dielectric cloaks \cite{ref_cloak}. Moreover, the GSM can be extracted with considerable flexibility---Besides tailored MoM and FEM implementations \cite{ref_myGSM,ref_3D_FEM,ref_GSM_PO}, the GSM can also be constructed using any solver capable of providing port-radiation and plane-wave scattering data. This paper provides a wrapper to demonstrate this capability through MATLAB--FEKO  co-simulation \cite{ref_github}. 

A second contribution of this work is the explicit derivation and implementation of the SSO formulation for a broad range of spherical media, including radially uniaxial anisotropic layers, multilayer configurations, and continuously varying radial profiles. Closed-form expressions are available for most homogeneous cases, while inhomogeneous profiles are handled by numerically integrating the governing radial equations using MATLAB packages such as \texttt{ODE45} and \texttt{ODE15s}. All scripts are openly released to support reproducibility and adaptation to other applications \cite{ref_github}. We note that the present SSO implementation is limited to radially uniaxial anisotropy; a general formulation for fully anisotropic, spatially varying spherical media remains an open challenge. This limitation concerns only SSO evaluation and does not affect the general GSM–SSO modeling paradigm, which remains broadly applicable to antennas embedded in spherical environments.

\begin{figure}[!t]
  \centering
  \includegraphics[]{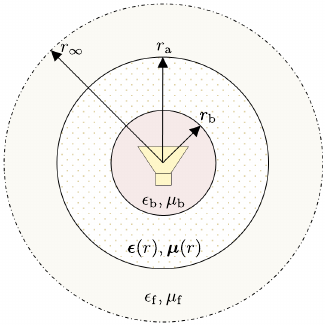}
  \caption{
    Canonical model for an embedded or implantable antenna. The antenna is enclosed within a spherical bubble of radius $r_\mathrm{b}$ and constitutive parameters $\epsilon_\mathrm{b}$ and $\mu_\mathrm{b}$. This bubble is surrounded by a spherical shell of outer radius $r_\mathrm{a}$, whose material may be radially inhomogeneous and uniaxially anisotropic, described by spatially varying tensors $\vec{\epsilon}(r)$ and $\vec{\mu}(r)$. The exterior region is homogeneous free space, characterized by $\epsilon_\mathrm{f}$ and $\mu_\mathrm{f}$, typically taken as 1.}
  \label{fGeneralScheme}
\end{figure}

\section{Hybrid Framework of Generalized Scattering Matrix and Spherical Scattering Operator}
\label{SecII}

We begin by considering the canonical configuration illustrated in Fig.~\ref{fGeneralScheme}. An antenna is enclosed within a spherical layer (hereafter referred to as the ``bubble'') of radius $r_\mrm{b}$. This bubble, in turn, is embedded inside a larger spherical region of radius $r_\mrm{a}$, which may consist of multiple layers with inhomogeneous and potentially anisotropic electromagnetic properties. The space beyond $r_\mrm{a}$ is free space. While both the bubble and the exterior region are assumed homogeneous and isotropic, the intermediate shell may exhibit arbitrary radial variations in its constitutive tensors $\vec{\epsilon}(r)$ and $\vec{\mu}(r)$. Such a geometry is widely encountered in practice, for example in bioelectromagnetic scenarios where an implanted antenna or sensor is enclosed by a protective coating and surrounded by biological tissues.

Inside the bubble, the electromagnetic fields are represented using regular and outgoing vector spherical wave functions (VSWFs), namely
\begin{equation}
  \label{eqn1}
  \begin{split}
     &\vec{E}\left( \vec{r} \right) =k_\mrm{b}\sqrt{Z_\mrm{b}}\sum_n{a_n^\mrm{b} \vec{u}_{n}^{\left( 1 \right)}\left( k_\mrm{b}\vec{r} \right) +f_n^\mrm{b}\vec{u}_{n}^{\left( 4 \right)}\left( k_\mrm{b}\vec{r} \right)}\\
     &\vec{H}\left( \vec{r} \right) =\frac{\j k_\mrm{b}}{\sqrt{Z_\mrm{b}}}\sum_n{a_n^\mrm{b}\vec{u}_{\bar n}^{\left( 1 \right)}\left( k_\mrm{b}\vec{r} \right) +f_n^\mrm{b}\vec{u}_{\bar n}^{\left( 4 \right)}\left( k_\mrm{b}\vec{r} \right)}
  \end{split}
\end{equation}
where $k_\mrm{b} = k_0 \sqrt{\mu_\mrm{b}\epsilon_\mrm{b}}$ and $Z_\mrm{b} = Z_0 \sqrt{\mu_\mrm{b}/\epsilon_\mrm{b}}$. $k_0$ and $Z_0$ are vacuum wave number and impedance, respectively. The VSWFs $\vec{u}_{n}^{\left( p \right)}$ are defined in Appendix~\ref{APP_A}, with superscripts $p=1$ and $p=4$ denoting regular and outgoing waves, respectively. The multi-index $n\to \tau\sigma m l$ includes polarization type $\tau =\left\{ 1,2 \right\}$ (TE, TM), parity $\sigma =\left\{ e,o \right\}$, degree $l=\left\{ 1,2,\cdots L_{\max} \right\}$, and order $m=\left\{ 0,1,\cdots ,l \right\}$. A bar over $n$ exchanges the TE and TM types \cite[Ch. 7]{ref_scattering_theory}.

A similar spherical-wave expansion holds in the outer free-space region, replacing $k_\mrm{b}$, $Z_\mrm{b}$ with $k_\mrm{f}$, $Z_\mrm{f}$:
\begin{equation}
  \begin{split}
     &\vec{E}\left( \vec{r} \right) =k_\mrm{f}\sqrt{Z_\mrm{f}}\sum_n{a_n^\mrm{f} \vec{u}_{n}^{\left( 1 \right)}\left( k_\mrm{f}\vec{r} \right) +f_n^\mrm{f}\vec{u}_{n}^{\left( 4 \right)}\left( k_\mrm{f}\vec{r} \right)}\\
     &\vec{H}\left( \vec{r} \right) =\frac{\j k_\mrm{f}}{\sqrt{Z_\mrm{f}}}\sum_n{a_n^\mrm{f}\vec{u}_{\bar n}^{\left( 1 \right)}\left( k_\mrm{f}\vec{r} \right) +f_n^\mrm{f}\vec{u}_{\bar n}^{\left( 4 \right)}\left( k_\mrm{f}\vec{r} \right)}
  \end{split}
\end{equation}

\begin{figure}[!t]
  \centering
  \includegraphics[]{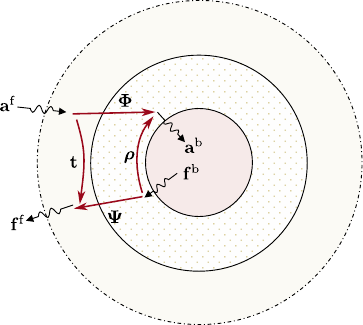}
  \caption{Schematic of wave propagation and interaction in a layered spherical system. Regular waves incident from free space ($\mathbf{a}^\mathrm{f}$) are partially reflected and partially transmitted into the inner bubble. Conversely, outgoing waves radiated by the antenna inside the bubble ($\mathbf{f}^\mathrm{b}$) experience partial transmission and reflection at the spherical shell. The operators $\mathbf{t}$, $\M{\Phi}$, $\M{\Psi}$, and $\vec{\rho}$ describe these propagation pathways}
  \label{fSpherical_operator}
\end{figure}

To relate the wave coefficients on both sides of the shell, we employ the bidirectional spherical scattering operator (SSO) \cite[Ch 3.5]{ref_InhomogeneousMedia}, \cite{ref_Capek1,ref_Capek2}. This operator characterizes the shell alone and is independent of the antenna, as demonstrated in Fig.~\ref{fSpherical_operator}:
\begin{equation}
  \label{eqn3}
   \begin{bmatrix}
	\mathbf{f}^{\mathrm{f}}\\
	\mathbf{a}^{\mathrm{b}}\\
\end{bmatrix} = \begin{bmatrix}
	\mathbf{t}&		\mathbf{\Psi }\\
	\mathbf{\Phi }&		\bm{\rho }\\
\end{bmatrix}  \begin{bmatrix}
	\mathbf{a}^{\mathrm{f}}\\
	\mathbf{f}^{\mathrm{b}}\\
\end{bmatrix} 
\end{equation}
Here, $\mat{t}$ denotes the free-space transition matrix (T-matrix); $\mat{\Phi}$ and $\mat{\Psi}$ describe inward and outward transmission; and $\bm{\rho}$ models internal reflections. In the absence of any radiating structure inside the bubble ($\mat{f}^\mathrm{b} = \M 0$), this relation naturally collapses to the standard external T-matrix formulation \cite{ref_Waterman}. 

\begin{figure}[!t]
  \centering
  \includegraphics[]{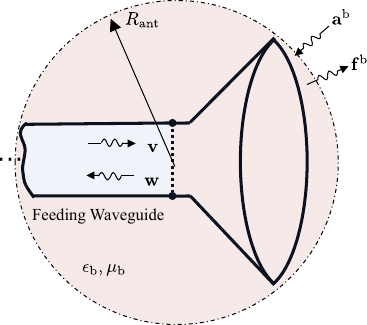}
  \caption{Antenna in a bubble-material filled free space. Regular incident and outgoing spherical waves in the bubble are $\mathbf{a}^\mathrm{b}$ and $\mathbf{f}^\mathrm{b}$, respectively. The antenna is interfaced with a waveguide supporting inward and outward modal amplitudes $\mathbf{v}$ and $\mathbf{w}$.}
  \label{fGSM}
  \vspace{-0.1in}
\end{figure}

In the present problem, the antenna is located inside the bubble. By the equivalence principle, the electromagnetic fields within the bubble can be treated as those of an antenna radiating in an unbounded medium identical to the bubble material. Accordingly, the spherical-wave behavior inside the bubble is governed by the antenna's free-space GSM. In this setting, the antenna is excited simultaneously by the guided-mode inputs $\mat{v}$ and the regular spherical waves $\mat{a}^\mathrm{b}$, and produces the waveguide response $\mat{w}$ as well as the outgoing spherical-wave coefficients $\mat{f}^\mathrm{b}$, as illustrated in Fig.~\ref{fGSM}. These quantities are related through the source–scattering formalism \cite[Ch.~2.3.5]{ref_sph_near_measure}:
\begin{equation}
  \label{eqn4}
  \begin{bmatrix}
	\mathbf{w}\\
	\mathbf{f}^{\mathrm{b}}\\
\end{bmatrix} = \begin{bmatrix}
	\mathbf{\Gamma }&		\frac{1}{2}\mathbf{R}\\
	\mathbf{T}&		\frac{1}{2}\left( \mathbf{S}-\mathbf{1} \right)\\
\end{bmatrix}  \begin{bmatrix}
	\mathbf{v}\\
	\mathbf{a}^{\mathrm{b}}\\
\end{bmatrix}
\end{equation}
where $\mat{\Gamma}$, $\mat{R}$, $\mat{T}$, and $\mat{S}$ are the port, reception, transmission, and scattering matrices, respectively.

In practical electromagnetic scenarios, the external excitation vector $\mathbf{a}^{\mathrm{f}}$ (\eg due to an incident plane wave) and the port excitation vector $\mathbf{v}$ are known. Thus, by combining \eqref{eqn3} and \eqref{eqn4} and eliminating the intermediate variables $\mathbf{a}^{\mathrm{b}}$ and $\mathbf{f}^{\mathrm{b}}$, the desired quantities $\mathbf{w}$ and $\mathbf{f}^{\mathrm{f}}$ can be expressed in closed form. To this end, rewriting \eqref{eqn4} yields
\begin{equation}
  \begin{cases}
	\mathbf{\Gamma v}+\frac{1}{2}\mathbf{Ra}^{\mathrm{b}}=\mathbf{w}\\
	\mathbf{Tv}+\frac{1}{2}\left( \mathbf{S}-\mathbf{1} \right) \mathbf{a}^{\mathrm{b}}=\mathbf{f}^{\mathrm{b}}\\
\end{cases}
\end{equation}
and substituting $\M a^\mrm{b} =  \mathbf{\Phi a}^{\mathrm{f}}+\bm{\rho}\M f^{\mathrm{b}}$ from \eqref{eqn3} leads to
\begin{equation}
  \label{eqn6}
  \begin{cases}
	\mathbf{\Gamma v}+\frac{1}{2}\mathbf{R}\left( \mathbf{\Phi a}^{\mathrm{f}}+\bm{\rho}\M f^{\mathrm{b}} \right) =\mathbf{w}\\
	\mathbf{Tv}+\frac{1}{2}\left( \mathbf{S}-\mathbf{1} \right) \left( \mathbf{\Phi a}^{\mathrm{f}}+\bm{\rho }\M f^{\mathrm{b}} \right) =\mathbf{f}^{\mathrm{b}}.\\
\end{cases}
\end{equation}
Solving the second equation for $\M f^\mrm{b}$ gives 
\begin{equation}
  \label{eq11}
  \mathbf{f}^{\mathrm{b}}=\mathbf{M}^{-1}\left[ \mathbf{Tv}+\frac{1}{2}\left( \mathbf{S}-\mathbf{1} \right) \mathbf{\Phi a}^{\mathrm{f}} \right] 
\end{equation}
where $\mathbf{M}=\left[ \mathbf{1}-\frac{1}{2}\left( \mathbf{S}-\mathbf{1} \right) \bm{\rho } \right]$. Substituting this into the expression for $\mat{w}$ in \eqref{eqn6} yields
\begin{equation}
  \label{eqn8}
  \begin{split}
    &\mathbf{w}=\left[ \mathbf{\Gamma }+\frac{1}{2}\mathbf{R}\bm{\rho} \M{M}^{-1}\mathbf{T} \right] \mathbf{v}\\
    &\quad+\frac{1}{2}\mathbf{R}\left[ \mathbf{\Phi }+\bm{\rho} \M{M}^{-1}\frac{1}{2}\left( \mathbf{S}-\mathbf{1} \right) \mathbf{\Phi } \right] \mathbf{a}^{\mathrm{f}}.
  \end{split}
\end{equation}

Using $\mathbf{f}^{\mathrm{f}} = \mathbf{t a}^{\mathrm{f}} + \mathbf{\Psi f}^{\mathrm{b}}$ and substituting \eqref{eq11} gives
\begin{equation}
  \label{eqn9}
  \mathbf{f}^{\mathrm{f}}=\mathbf{\Psi M}^{-1}\mathbf{Tv}+\left[ \mathbf{t}+\mathbf{\Psi M}^{-1}\frac{1}{2}\left( \mathbf{S}-\mathbf{1} \right) \mathbf{\Phi } \right] \mathbf{a}^{\mathrm{f}}.
\end{equation}

Finally, \eqref{eqn8}--\eqref{eqn9} can be assembled as
\begin{equation}
  \label{eqn10}
   \begin{bmatrix}
	\mathbf{w}\\
	\mathbf{f}^{\mathrm{f}}\\
\end{bmatrix} = \begin{bmatrix}
	\tM{\Gamma}&		\frac{1}{2}\tM{R}\\
	\tM{T}&		\frac{1}{2}\left( \tM{S}-\mathbf{1} \right)\\
\end{bmatrix} \begin{bmatrix}
	\mathbf{v}\\
	\mathbf{a}^{\mathrm{f}}\\
\end{bmatrix}
\end{equation}
where
\begin{equation}
  \begin{split}
    &\tM{\Gamma}=\mathbf{\Gamma }+\frac{1}{2}\mathbf{R}\bm{\rho} \M{M}^{-1}\mathbf{T}\\
    &\tM{R}=\mathbf{R}\left[ \mathbf{\Phi }+\bm{\rho} \M{M}^{-1}\frac{1}{2}\left( \mathbf{S}-\mathbf{1} \right) \mathbf{\Phi } \right] \\
    &\tM{T}=\mathbf{\Psi M}^{-1}\mathbf{T}\\
    &\tM{S}=\mathbf{1}+2\mathbf{t}+\mathbf{\Psi M}^{-1}\left( \mathbf{S}-\mathbf{1} \right) \mathbf{\Phi }.
  \end{split}
\end{equation}
These matrix blocks inherit the physical meaning of their free-space counterparts in \eqref{eqn4}, while embedding the full electromagnetic response of the surrounding spherical medium. They therefore form the GSM of the full antenna--stratified-sphere configuration.

Accuracy is ensured by choosing the spherical-harmonic truncation order as
\begin{equation*}
  L_{\max}=\lceil \varrho +2\sqrt[3]{\varrho }+3 \rceil ,\quad
   \varrho = \max\{k_{\mathrm{f}} r_{\mathrm{a}},\, k_{\mathrm{b}} R_{\mathrm{ant}}\}
\end{equation*}
following~\cite{ref_Sph_deg_trunction}. When $k_{\mathrm{b}} R_{\mathrm{ant}} < k_{\mathrm{f}} r_{\mathrm{a}}$, the blocks $\mathbf{R}$, $\mathbf{T}$, and $\mathbf{S}-\mathbf{1}$ can be zero-padded to match the SSO matrix dimensions, which holds in the vast majority of practical situations and reduces computational cost for determining the GSM.

\section{Determination of the Scattering Operators}
\label{Sec_III}

To derive the spherical scattering operators $\mathbf{t},\mathbf{\Psi},\mathbf{\Phi},\boldsymbol{\rho}$, we first expand the electromagnetic fields within the intermediate region in terms of spherical vector basis functions ${\vec w_n(\vec r)}$,
\begin{equation}
  \label{eqn12}
  \begin{split}
    &\vec{E}\left( \vec{r} \right) =\sum_n{\alpha _n\vec{w}_n\left( \vec{r} \right)}\\
	  &\vec{H}\left( \vec{r} \right) =\frac{\j}{k_0Z_0}\sum_n{\alpha _n\vec{\mu }^{-1}\left( \vec{r} \right) \cdot \left( \nabla \times \vec{w}_n\left( \vec{r} \right) \right)}
  \end{split}
\end{equation}
where each basis function is expressed as a radial combination of vector spherical harmonics,
\begin{equation}
  \vec{w}_n=I_{1n}\left( r \right) \vec{A}_{1n}+I_{2n}\left( r \right) \vec{A}_{2n}+I_{3n}\left( r \right) \vec{A}_{3n}
\end{equation}
with $\vec A_{\tau n}$ defined in Appendix~\ref{APP_A}. At this stage, no assumption has yet been made regarding the material tensors $\boldsymbol{\epsilon}(r)$ and $\boldsymbol{\mu}(r)$; the expansion remains completely general.

To obtain explicit and computationally tractable expressions for the radial functions $I_{1n}(r)$, $I_{2n}(r)$, and $I_{3n}(r)$, we next specialize to a physically important class of materials: radially uniaxial media with arbitrary radial inhomogeneity. In this case, the constitutive tensors take the form
\begin{equation}
  \label{eqn14}
  \begin{split}
    &\vec{\epsilon }\left( r \right) =\epsilon _{\bot}\left( r \right) \left( \mathbf{1}_3-\hat{\vec{r}}\hat{\vec{r}} \right) +\epsilon _r\left( r \right) \hat{\vec{r}}\hat{\vec{r}}\\
    &\vec{\mu }\left( r \right) =\mu _{\bot}\left( r \right) \left( \mathbf{1}_3-\hat{\vec{r}}\hat{\vec{r}} \right) +\mu _r\left( r \right) \hat{\vec{r}}\hat{\vec{r}}.
  \end{split}
\end{equation}
where $\mat 1_3$ denotes the 3D identity dyadic. Under this symmetry, the TE ($\tau=1$) and TM ($\tau=2$) spherical modes are uncoupled, and the basis functions naturally separate into the canonical forms \cite[Ch.~8.4]{ref_scattering_theory}:
\begin{equation}
  \label{eqn15}
  \begin{split}
    &\vec{w}_{1n}=\frac{g\left( r \right)}{k_0r}\vec{A}_{1n}\left( \hat{\vec{r}} \right)\\
    &\vec{w}_{2n}=\frac{h^{\prime}\left( r \right)}{k_{0}^{2}r\epsilon _{\bot}\left( r \right)}\vec{A}_{2n}\left( \hat{\vec{r}} \right) +\frac{h\left( r \right) \sqrt{l\left( l+1 \right)}}{k_{0}^{2}\epsilon _r\left( r \right) r^2}\vec{A}_{3n}\left( \hat{\vec{r}} \right)
  \end{split}
\end{equation}
with radial functions governed by the second-order ODE system
\begin{equation}
  \label{eqn16}
  \begin{cases}
	g^{\prime \prime}\left( r \right) +p_1\left( r \right) g^{\prime}\left( r \right) +q_1\left( r \right) g\left( r \right) =0\\
	h^{\prime \prime}\left( r \right) +p_2\left( r \right) h^{\prime}\left( r \right) +q_2\left( r \right) h\left( r \right) =0\\
  \end{cases}
\end{equation}
where
\begin{equation}
  \begin{split}
    &p_1=-\frac{\mu _{\bot}^{\prime}\left( r \right)}{\mu _{\bot}\left( r \right)},q_1=k_{0}^{2}\mu _{\bot}\left( r \right) \epsilon _{\bot}\left( r \right) -\frac{\mu _{\bot}\left( r \right) l\left( l+1 \right)}{\mu _r\left( r \right) r^2}\\
    &p_2=-\frac{\epsilon _{\bot}^{\prime}\left( r \right)}{\epsilon _{\bot}\left( r \right)},q_2=k_{0}^{2}\mu _{\bot}\left( r \right) \epsilon _{\bot}\left( r \right) -\frac{\epsilon _{\bot}\left( r \right) l\left( l+1 \right)}{\epsilon _r\left( r \right) r^2}.\quad
  \end{split}
\end{equation}

For homogeneous media, $\epsilon_r$, $\epsilon_\perp$, $\mu_r$, and $\mu_\perp$ are constant, so $p_1=p_2=0$. In this case, the radial functions admit closed-form Riccati-Bessel solutions
\begin{equation}
  \label{eqn18}
  \begin{split}
      &g_l\left( r \right)   =A\zeta _{L_1}\left( k_{\bot}r \right) +B\xi _{L_1}\left( k_{\bot}r \right)\\
	    &h_l\left( r \right)  =C\zeta _{L_2}\left( k_{\bot}r \right) +D\xi _{L_2}\left( k_{\bot}r \right)
  \end{split}
\end{equation}
where
\begin{equation}
  \begin{split}
	L_1&=\sqrt{\frac{\mu _{\bot}}{\mu _r}l\left( l+1 \right) +\frac{1}{4}}-\frac{1}{2}\\
	L_2&=\sqrt{\frac{\epsilon _{\bot}}{\epsilon _r}l\left( l+1 \right) +\frac{1}{4}}-\frac{1}{2}
\end{split}
\end{equation}
$k_\perp=k_0\sqrt{\epsilon_\perp\mu_\perp}$ and the coefficients $A,B,C,D$ are determined uniquely from boundary continuity conditions. The isotropic case follows as a special case with $L_1=L_2=l$.

For general radially varying media, analytic solutions of \eqref{eqn16} are not available; however, numerical integration remains efficient. In the examples presented in this work, the radial profiles can be clearly identified as non-stiff, and the ODEs are therefore solved using MATLAB's \texttt{ODE45} with a recommended \texttt{RelTol} of $10^{-6}$. For problems in which users cannot readily determine the stiffness a priori, a practical strategy is to first solve the ODE with \texttt{ODE45} and then perform a cross-check using the stiff solver \texttt{ODE15s}. Additional confidence can be obtained by adjusting the tolerance---\eg repeating the highest spherical-harmonic order with two different \texttt{RelTol} values---to verify convergence. 

If the medium deviates from the radially uniaxial form in \eqref{eqn14}, TE--TM decoupling no longer holds and the separation in \eqref{eqn15} breaks down. Closed-form spherical solutions are then generally unavailable, except for the restrictive case of spatially invariant Cartesian-aligned anisotropy, and even then the resulting expressions remain cumbersome \cite{ref_Stout_Mie,ref_Rafaei_DGF}. For inhomogeneous Cartesian-aligned anisotropy, one could in principle formulate a coupled ODE system analogous to \eqref{eqn16}; however, the breakdown of variable separation presents a fundamental algebraic obstruction for this entire class of problems. Consequently, and to the best of our knowledge, no tractable general formulation has been reported, leaving this avenue open for future investigation. Importantly, these complications affect only the evaluation of $(\mathbf{t},\mathbf{\Psi},\mathbf{\Phi},\boldsymbol{\rho})$; the GSM-SSO framework itself remains unchanged.

\subsection{Determination of \texorpdfstring{$\M t$ matrix and $\M\Phi$ matrix}{}}
\label{SecIII_A}

We first consider the case where the bubble contains no internal source, \ie $\M f^\mathrm{b} = \M 0$. Under this condition, \eqref{eqn3} reduces to $\M f^\mathrm{f} = \M t \M a^\mathrm{f}$ and $\M a^\mathrm{b} = \M \Phi \M a^\mathrm{f}$.

At the inner interface $r = r_\mathrm{b}$ (with $x_\mathrm{b}\triangleq k_\mathrm{b}r_\mathrm{b}$), continuity of tangential fields (\ie $\hat{\vec r}\times$\eqref{eqn1} and $\hat{\vec r}\times$\eqref{eqn12}), together with the orthogonality of the vector spherical harmonics $\vec A_{\tau n}$, yields the following boundary relations for TE ($\tau=1$) modes \cite[Ch.~8.4]{ref_scattering_theory}:
\begin{equation}
  \label{eqn20}
  \begin{split}
	&k_{\mathrm{b}}\sqrt{Z_{\mathrm{b}}}a_{1n}^{\mathrm{b}}R_{1l}^{\left( 1 \right)}\left( x_\mrm{b} \right) =\alpha _{1n}\frac{g\left( r_{\mathrm{b}} \right)}{k_0r_{\mathrm{b}}}\\
	&\frac{k_{\mathrm{b}}}{\sqrt{Z_{\mathrm{b}}}}a_{1n}^{\mathrm{b}}R_{2l}^{\left( 1 \right)}\left( x_\mrm{b} \right) =\alpha _{1n}\frac{1}{Z_0}\frac{g^{\prime}\left( r_{\mathrm{b}} \right)}{\mu _{\bot}\left( r_{\mathrm{b}} \right) k_{0}^{2}r_{\mathrm{b}}}
\end{split}
\end{equation}
and for TM ($\tau=2$) modes: 
\begin{equation}
  \label{eqn21}
  \begin{split}
	&k_{\mathrm{b}}\sqrt{Z_{\mathrm{b}}}a_{2n}^{\mathrm{b}}R_{2l}^{\left( 1 \right)}\left( x_\mrm{b} \right) =\alpha _{2n}\frac{h^{\prime}\left( r_{\mathrm{b}} \right)}{k_{0}^{2}r_{\mathrm{b}}\epsilon _{\bot}\left( r_{\mathrm{b}} \right)}\\
	&\frac{k_{\mathrm{b}}}{\sqrt{Z_{\mathrm{b}}}}a_{2n}^{\mathrm{b}}R_{1l}^{\left( 1 \right)}\left( x_\mrm{b}\right) =\alpha _{2n}\frac{1}{Z_0}\frac{h\left( r_{\mathrm{b}} \right)}{k_0r_{\mathrm{b}}}.
\end{split} 
\end{equation}
Here, $n\to \sigma m l$, and thus the relations are independent of $\sigma m$ indices. From \eqref{eqn20}--\eqref{eqn21}, one obtains the initial conditions
\begin{equation}
  \label{eqn22}
  \begin{split}
    &\frac{g^{\prime}\left( r_{\mathrm{b}} \right)}{g\left( r_{\mathrm{b}} \right)}=\frac{k_{\mathrm{b}} \mu _{\bot}\left( r_{\mathrm{b}} \right) }{\mu _{\mathrm{b}}}\frac{\zeta _{l}^{\prime}\left( x_\mrm{b} \right)}{\zeta _l\left( x_\mrm{b} \right)}\\
    &\frac{h^{\prime}\left( r_{\mathrm{b}} \right)}{h\left( r_{\mathrm{b}} \right)}=\frac{k_{\mathrm{b}}\epsilon _{\bot}\left( r_{\mathrm{b}} \right)}{\epsilon _{\mathrm{b}}}\frac{\zeta _{l}^{\prime}\left( x_\mrm{b} \right)}{\zeta _l\left( x_\mrm{b} \right)}
  \end{split}
\end{equation}
where $g(r_\mathrm{b}) = h(r_\mathrm{b}) = 1$ sets the normalization. These conditions uniquely determine $g(r)$, $h(r)$ and their derivatives in \eqref{eqn16} (in the analytical case, they fix the undetermined coefficient $A,B,C,D$ in \eqref{eqn18}). The relations between $R_{\tau l}^{(p)}$ and Riccati-Bessel functions $\zeta_l$, $\xi_l$ appear in Appendix~\ref{APP_A}.

At the outer interface $r = r_\mathrm{a}$, enforcing continuity of tangential fields and noting $f_{\tau n}^\mathrm{f}=t_{\tau n}a_{\tau n}^\mathrm{f}$ yields
\begin{equation}
\label{eqn23}
\begin{aligned}
&\alpha _{1n}\frac{g(r_{\mathrm{a}})}{k_0 r_{\mathrm{a}}}
= k_{\mathrm{f}}\sqrt{Z_{\mathrm{f}}}\, a_{1n}^{\mathrm{f}}\left[
R_{1l}^{(1)}(x_{\mathrm{a}})+t_{1n} R_{1l}^{(4)}(x_{\mathrm{a}})
\right]\quad \\[6pt]
&\alpha _{1n}\frac{1}{Z_0}\frac{g'(r_{\mathrm{a}})}{\mu _{\bot}(r_{\mathrm{a}}) k_{0}^{2} r_{\mathrm{a}}}
= \frac{k_{\mathrm{f}} a_{1n}^{\mathrm{f}}}{\sqrt{Z_{\mathrm{f}}}}\left[
R_{2l}^{(1)}(x_{\mathrm{a}})+t_{1n}R_{2l}^{(4)}(x_{\mathrm{a}})
\right]\quad
\end{aligned}
\end{equation}
and 
\begin{equation}
  \label{eqn24}
  \begin{split}
	&\alpha _{2n}\frac{h^{\prime}\left( r_{\mathrm{a}} \right)}{k_{0}^{2}r_{\mathrm{a}}\epsilon _{\bot}\left( r_{\mathrm{a}} \right)}=k_{\mathrm{f}}\sqrt{Z_{\mathrm{f}}}a_{2n}^{\mathrm{f}}\left[ R_{2l}^{\left( 1 \right)}\left( x_{\mathrm{a}} \right) +t_{2n}R_{2l}^{\left( 4 \right)}\left( x_{\mathrm{a}} \right) \right]\\
	&\alpha _{2n}\frac{1}{Z_0}\frac{h\left( r_{\mathrm{a}} \right)}{k_0r_{\mathrm{a}}}=\frac{k_{\mathrm{f}}a_{2n}^{\mathrm{f}}}{\sqrt{Z_{\mathrm{f}}}}\left[ R_{1l}^{\left( 1 \right)}\left( x_{\mathrm{a}} \right) +t_{2n}R_{1l}^{\left( 4 \right)}\left( x_{\mathrm{a}} \right) \right]
\end{split}
\end{equation}
with $x_\mathrm{a}\triangleq k_\mathrm{f}r_\mathrm{a}$.

By dividing the two rows in \eqref{eqn23}, the coefficients $\alpha_{1n}$ and $a^\mrm{f}_{1n}$ cancel, yielding a single equation dependent only on $t_{1n}$:
\begin{equation}
  \label{eqn25}
  \frac{k_0Z_0\mu _{\bot}(r_{\mathrm{a}})g(r_{\mathrm{a}})}{Z_{\mathrm{f}}g^\prime (r_{\mathrm{a}})}=\frac{ R_{1l}^{(1)}(x_{\mathrm{a}})+t_{1n}R_{1l}^{(4)}(x_{\mathrm{a}})}{ R_{2l}^{(1)}(x_{\mathrm{a}})+t_{1n}R_{2l}^{(4)}(x_{\mathrm{a}})}.
\end{equation}
Repeating the same procedure on \eqref{eqn24} gives
\begin{equation}
  \label{eqn26}
  \frac{Z_0h^{\prime}\left( r_{\mathrm{a}} \right)}{k_0Z_{\mathrm{f}}\epsilon _{\bot}\left( r_{\mathrm{a}} \right) h\left( r_{\mathrm{a}} \right)}=\frac{R_{2l}^{\left( 1 \right)}\left( x_{\mathrm{a}} \right) +t_{2n}R_{2l}^{\left( 4 \right)}\left( x_{\mathrm{a}} \right)}{R_{1l}^{\left( 1 \right)}\left( x_{\mathrm{a}} \right) +t_{2n}R_{1l}^{\left( 4 \right)}\left( x_{\mathrm{a}} \right)}
\end{equation}

Equations \eqref{eqn25} and \eqref{eqn26} then lead to explicit semi-closed forms for the transition coefficients $t_{1n}$ and $t_{2n}$, respectively. Using $\frac{k_0 Z_0}{Z_\mrm{f}}=\frac{k_\mrm{f}}{\mu_\mrm{f}}$ and $\frac{Z_0}{k_0 Z_\mrm{f}}=\frac{\epsilon_\mrm{f}}{k_\mrm{f}}$, we arrive at
\begin{equation}
  \label{eqn27}
  \begin{split}
    &t_{1n}=-\frac{\displaystyle  \frac{\mu _{\mathrm{f}}}{k_{\mathrm{f}}\mu _{\bot}\left( r_{\mathrm{a}} \right)}\frac{g^{\prime}\left( r_{\mathrm{a}} \right)}{g\left( r_{\mathrm{a}} \right)}\zeta _l\left( x_\mrm{a} \right) -\zeta _{l}^{\prime}\left( x_\mrm{a} \right)}{\displaystyle  \frac{\mu _{\mathrm{f}}}{k_{\mathrm{f}}\mu _{\bot}\left( r_{\mathrm{a}} \right)}\frac{g^{\prime}\left( r_{\mathrm{a}} \right)}{g\left( r_{\mathrm{a}} \right)}\xi _l\left( x_\mrm{a} \right) -\xi _{l}^{\prime}\left( x_\mrm{a} \right)}\\
    &t_{2n}=-\frac{\displaystyle  \frac{\epsilon _{\mathrm{f}}}{k_{\mathrm{f}}\epsilon _{\bot}\left( r_{\mathrm{a}} \right)}\frac{h^{\prime}\left( r_{\mathrm{a}} \right)}{h\left( r_{\mathrm{a}} \right)}\zeta _l\left( x_\mrm{a} \right) -\zeta _{l}^{\prime}\left( x_\mrm{a} \right)}{\displaystyle  \frac{\epsilon _{\mathrm{f}}}{k_{\mathrm{f}}\epsilon _{\bot}\left( r_{\mathrm{a}} \right)}\frac{h^{\prime}\left( r_{\mathrm{a}} \right)}{h\left( r_{\mathrm{a}} \right)}\xi _l\left( x_\mrm{a} \right) -\xi _{l}^{\prime}\left( x_\mrm{a} \right)}
  \end{split}
\end{equation}

To obtain $\Phi_{1n}$ and $\Phi_{2n}$, we form the ratios $\Phi_{1n}=a_{1n}^\mrm{b}/a_{1n}^\mrm{f}$ and $\Phi_{2n}=a_{2n}^\mrm{b}/a_{2n}^\mrm{f}$. Eliminating $\alpha_{1n}$ from the first rows of \eqref{eqn20} and \eqref{eqn23}, and eliminating $\alpha_{2n}$ from the first rows of \eqref{eqn21} and \eqref{eqn24}, yields
\begin{equation}
  \label{eqn28}
  \begin{split}
    &\Phi _{1n}=\frac{a_{1n}^{\mathrm{b}}}{a_{1n}^{\mathrm{f}}}=\frac{g\left( r_{\mathrm{b}} \right)}{g\left( r_{\mathrm{a}} \right)}\frac{\sqrt{Z_{\mathrm{f}}}}{\sqrt{Z_{\mathrm{b}}}}\frac{\zeta _l\left( x_\mrm{a} \right) +t_{1n}\xi _l\left( x_\mrm{a} \right)}{\zeta _l\left( x_\mrm{b} \right)}\\
    &\Phi _{2n}=\frac{a_{2n}^{\mathrm{b}}}{a_{2n}^{\mathrm{f}}}=\frac{h^{\prime}\left( r_{\mathrm{b}} \right)}{h^{\prime}\left( r_{\mathrm{a}} \right)}\frac{\epsilon _{\bot}\left( r_{\mathrm{a}} \right)}{\epsilon _{\bot}\left( r_{\mathrm{b}} \right)}\frac{\sqrt{Z_{\mathrm{f}}}}{\sqrt{Z_{\mathrm{b}}}}\frac{\zeta _{l}^{\prime}\left( x_\mrm{a} \right) +t_{2n}\xi _{l}^{\prime}\left( x_\mrm{a} \right)}{\zeta _{l}^{\prime}\left( x_\mrm{b} \right)}.\quad
  \end{split}
\end{equation}

Therefore, as seen from \eqref{eqn27} and \eqref{eqn28}, both $\mathbf{t}$ and $\mathbf{\Phi}$ are diagonal in $\tau l$ indices.

\subsection{Determination of \texorpdfstring{$\bm \rho$ matrix and $\M \Psi$ matrix}{}}
\label{SecIII_B}

The derivation of the matrices $\bm{\rho}$ and $\mathbf{\Psi}$ follows the same procedure, now assuming no external excitation, \ie $\mathbf{a}^\mathrm{f} = \mathbf{0}$. Under this condition, $\mathbf{a}^\mathrm{b} = \bm{\rho} \mathbf{f}^\mathrm{b}$ and $\mathbf{f}^\mathrm{f} = \mathbf{\Psi} \mathbf{f}^\mathrm{b}$. The initial conditions are imposed at the outer interface $r = r_\mathrm{a}$:
\begin{equation}
  \label{eqn29}
  \begin{split}
    &\frac{g^{\prime}\left( r_{\mathrm{a}} \right)}{g\left( r_{\mathrm{a}} \right)}=\frac{k_{\mathrm{f}}\mu _{\bot}\left( r_{\mathrm{a}} \right)}{\mu _{\mathrm{f}}}\frac{\xi _{l}^{\prime}\left( k_{\mathrm{f}}r_{\mathrm{a}} \right)}{\xi _l\left( k_{\mathrm{f}}r_{\mathrm{a}} \right)}\\
    &\frac{h^{\prime}\left( r_{\mathrm{a}} \right)}{h\left( r_{\mathrm{a}} \right)}=\frac{k_{\mathrm{f}}\epsilon _{\bot}\left( r_{\mathrm{a}} \right)}{\epsilon _{\mathrm{f}}}\frac{\xi _{l}^{\prime}\left( k_{\mathrm{f}}r_{\mathrm{a}} \right)}{\xi _l\left( k_{\mathrm{f}}r_{\mathrm{a}} \right)}
  \end{split}
\end{equation}
with the normalization $g(r_\mrm{a})=h(r_\mrm{a})=1$. Solving the ODE system \eqref{eqn16} with these conditions yields the radial functions $g(r)$ and $h(r)$, from which $\bm{\rho}$ and $\mathbf{\Psi}$ follow by the same algebraic steps as in the previous subsection.

It is worth noting that \eqref{eqn29} is the dual counterpart of \eqref{eqn22}, obtained through the mapping $r_\mrm{b}\leftrightarrow r_\mrm{a}$, $k_\mrm{b}\leftrightarrow k_\mrm{f}$, $\mu_\mrm{b}\leftrightarrow \mu_\mrm{f}$, $\epsilon_\mrm{b}\leftrightarrow\epsilon_\mrm{f}$, and $\zeta_l \leftrightarrow \xi_l$. Accordingly, the semi-closed-form expressions for $\rho_{\tau n}$ can be obtained by applying this duality to \eqref{eqn27}:
\begin{equation}
  \label{eqn30}
  \begin{split}
    &\rho _{1n}=-\frac{\displaystyle  \frac{\mu _{\mathrm{b}}}{k_{\mathrm{b}}\mu _{\bot}\left( r_{\mathrm{b}} \right)}\frac{g^{\prime}\left( r_{\mathrm{b}} \right)}{g\left( r_{\mathrm{b}} \right)}\xi _l\left( x_\mrm{b} \right) -\xi _{l}^{\prime}\left( x_\mrm{b} \right)}{\displaystyle  \frac{\mu _{\mathrm{b}}}{k_{\mathrm{b}}\mu _{\bot}\left( r_{\mathrm{b}} \right)}\frac{g^{\prime}\left( r_{\mathrm{b}} \right)}{g\left( r_{\mathrm{b}} \right)}\zeta _l\left( x_\mrm{b} \right) -\zeta _{l}^{\prime}\left( x_\mrm{b} \right)}\\
    &\rho _{2n}=-\frac{\displaystyle  \frac{\epsilon _{\mathrm{b}}}{k_{\mathrm{b}}\epsilon _{\bot}\left( r_{\mathrm{b}} \right)}\frac{h^{\prime}\left( r_{\mathrm{b}} \right)}{h\left( r_{\mathrm{b}} \right)}\xi _l\left( x_\mrm{b} \right) -\xi _{l}^{\prime}\left( x_\mrm{b} \right)}{\displaystyle  \frac{\epsilon _{\mathrm{b}}}{k_{\mathrm{b}}\epsilon _{\bot}\left( r_{\mathrm{b}} \right)}\frac{h^{\prime}\left( r_{\mathrm{b}} \right)}{h\left( r_{\mathrm{b}} \right)}\zeta _l\left( x_\mrm{b} \right) -\zeta _{l}^{\prime}\left( x_\mrm{b} \right)}.
  \end{split}
\end{equation}

Similarly, applying the same duality to \eqref{eqn28} (with the additional substitution $t_{\tau n}\leftrightarrow\rho_{\tau n}$) yields
\begin{equation}
  \label{eqn31}
  \begin{split}
    &\Psi _{1n}=\frac{g\left( r_{\mathrm{a}} \right)}{g\left( r_{\mathrm{b}} \right)}\frac{\sqrt{Z_{\mathrm{b}}}}{\sqrt{Z_{\mathrm{f}}}}\frac{\rho _{1n}\zeta _l\left( x_\mrm{b} \right) +\xi _l\left( x_\mrm{b} \right)}{\xi _l\left( x_\mrm{a} \right)}\\
    &\Psi _{2n}=\frac{h^{\prime}\left( r_{\mathrm{a}} \right)}{h^{\prime}\left( r_{\mathrm{b}} \right)}\frac{\epsilon _{\bot}\left( r_{\mathrm{b}} \right) \sqrt{Z_{\mathrm{b}}}}{\epsilon _{\bot}\left( r_{\mathrm{a}} \right) \sqrt{Z_{\mathrm{f}}}}\frac{\rho _{2n}\zeta _{l}^{\prime}\left( x_\mrm{b} \right) +\xi _{l}^{\prime}\left( x_\mrm{b} \right)}{\xi _{l}^{\prime}\left( x_\mrm{a} \right)}
  \end{split}
\end{equation}

When the bubble medium equals the surrounding free space, the reciprocity symmetry gives $\mathbf{\Phi} = \mathbf{\Psi}^t$.

\section{Numerical Validations}
\label{Sec_IV}

This section presents several representative examples to verify the accuracy, flexibility, and computational efficiency of the proposed method. The radiating structure under test is a multimode horn antenna identical to that in \cite{ref_myGSM}. The horn supports five waveguide modes over the 3.2--3.8~GHz band, namely $\mathrm{TE}_{10}$, $\mathrm{TE}_{20}$, $\mathrm{TE}_{01}$, $\mathrm{TE}_{11}$, and $\mathrm{TM}_{11}$. 

The antenna is embedded inside a stratified spherical medium, as illustrated in Fig.~\ref{fGeneralScheme}. Both the inner bubble and the outer region are vacuum (air); material inhomogeneity and anisotropy are confined to the intermediate spherical shell.

In all validation cases, the antenna GSM blocks $\mathbf{\Gamma}$, $\mathbf{R}$, $\mathbf{T}$, and $\mathbf{S}$ are obtained using a MoM-based procedure \cite{ref_myGSM} with 5901 RWG basis functions \cite{ref_RWG}.

\subsection{Homogeneous Isotropic Medium}

The first scenario considers a homogeneous, isotropic, lossy medium in the intermediate layer, with relative permittivity $\epsilon _r=\epsilon _{\bot}=5-0.5\j$ and permeability $\mu _r=\mu _{\bot}=1$. The inner and outer interface radii are $r_\mrm{b}=\mrm{150~mm}$ and $r_\mrm{a}=\mrm{180~mm}$, respectively.

For this special case, the radial ODEs in \eqref{eqn16} admit the closed-form Riccati-Bessel solutions in \eqref{eqn18}, and the transition matrices $\mathbf{t}$, $\mathbf{\Phi}$, $\boldsymbol{\rho}$, and $\mathbf{\Psi}$ follow analytically from \eqref{eqn27}--\eqref{eqn31}. The complete electromagnetic fields are subsequently assembled using \eqref{eqn10}.

\begin{figure}[!t]
  \centering
  \includegraphics[]{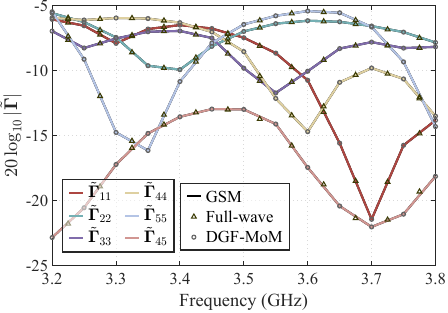}
  \caption{Port S-parameters of the multimode horn antenna embedded in a homogeneous isotropic spherical shell ($\epsilon_r=\epsilon_{\bot}=5-0.5\j$, $\mu_r=\mu_{\bot}=1$). Solid curves: proposed formulation. Circles: MoM with closed-form DGF. Modes 1--5 correspond to $\mathrm{TE}_{10}$, $\mathrm{TE}_{20}$, $\mathrm{TE}_{01}$, $\mathrm{TE}_{11}$, and $\mathrm{TM}_{11}$, respectively.}
\label{fSpara_Horn_singleIsotropic}
\end{figure}

Figure~\ref{fSpara_Horn_singleIsotropic} shows the S-parameters at the horn port. Curves labeled ``GSM'' are produced by the proposed formulation, while ``DGF-MoM'' corresponds to a method-of-moments implementation equipped with a closed-form spherical Green's function for this geometry. Both sets exhibit excellent agreement with full-wave FEKO simulations.

Figure~\ref{fSpara_Horn_singleIsotropic} plot the horn port S-parameters. The curves labeled ``GSM'' are produced using the proposed method, while those labeled ``DGF-MoM'' come from an MoM implementation equipped with a closed-form spherical Green's function tailored to this geometry. Both sets show excellent agreement with full-wave FEKO simulations.

Another quantity of interest is the far-field pattern, defined as
\begin{equation}
  \vec{F}\left( \hat{\vec{r}} \right) =\lim_{r\rightarrow \infty} re^{\j k_\mrm{f} r}\vec{E}\left( \hat{\vec{r}} \right) .
\end{equation}
For a field represented through the outgoing coefficients $\M f^\mrm{f}$, $\vec{F}\left( \hat{\vec{r}} \right)$, the far field follows from the asymptotic form of spherical vector waves \cite{ref_Capek1}:
\begin{equation}
  \label{eqn33}
  \vec{F}\left( \hat{\vec{r}} \right) =\sqrt{Z_{\mathrm{f}}}\sum_{n}{f_{n}^{\mathrm{f}}\mathrm{j}^{l+2-\tau}\vec{A}_{n}\left( \hat{\vec{r}} \right) }.
\end{equation}

\begin{figure}[!t]
  \centering
  \includegraphics[]{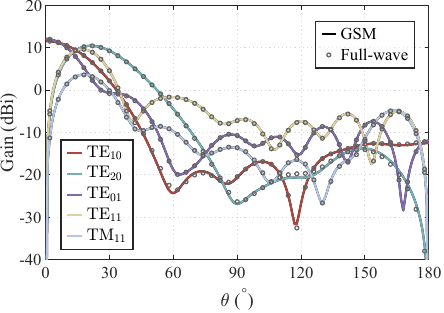}
  \caption{Far-field gain patterns on the $xoz$ plane at 3.2~GHz for the multimode horn antenna embedded in an isotropic homogeneous spherical shell. Each curve corresponds to a distinct waveguide excitation mode among $\mathrm{TE}_{10}$, $\mathrm{TE}_{20}$, $\mathrm{TE}_{01}$, $\mathrm{TE}_{11}$, and $\mathrm{TM}_{11}$.}
  \label{fGain_Horn_singleIsotropic}
\end{figure}

To compute the radiation patterns under waveguide excitation, we set $\mathbf{a}^\mathrm{f}=0$, which yields $\mathbf{f}^\mathrm{f}=\mathbf{T}\mathbf{v}$. Fig.~\ref{fGain_Horn_singleIsotropic} plots the gain in the $xoz$-plane at 3.2~GHz for each excitation mode (the $\mathrm{TE}_{01}$ pattern corresponds to the E-plane). Again, the GSM-SSO results closely match FEKO.

\begin{figure}[!t]
  \centering
  \includegraphics[]{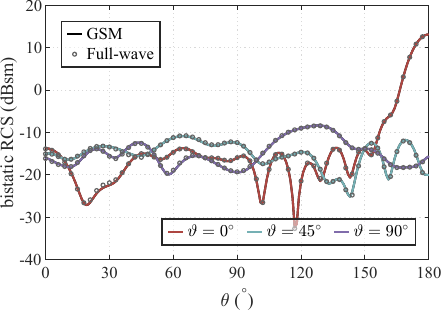}
  \caption{Bistatic RCS in the $xoz$-plane at 3.5~GHz for the multimode horn embedded in a homogeneous isotropic spherical shell. Results shown for plane-wave incidence at $\vartheta=0^\circ$, $45^\circ$, and $90^\circ$.}
  \label{fRCS_Horn_singleIsotropic}
\end{figure}

For scattering analysis, $\mathbf{v}=0$, leading to $\mathbf{f}^\mathrm{f}=\tfrac{1}{2}(\mathbf{S}-\mathbf{1})\mathbf{a}^\mathrm{f}$. The plane-wave excitation coefficients $\mathbf{a}^{\mathrm{f}}$ follow \cite[Appendix~A]{ref_myGSM}. Fig.~\ref{fRCS_Horn_singleIsotropic} shows the bistatic RCS for several incidence angles. The agreement with FEKO remains excellent, with minor discrepancies attributable to the discretization of the spherical shell.

\subsection{Homogeneous Anisotropic Medium}

We next examine a homogeneous uniaxially anisotropic shell with $\epsilon_r=2$, $\epsilon_{\bot}=5$, $\mu_r=1$, and $\mu_{\bot}=3$, using the same radii as before. For this medium, \eqref{eqn18} remains valid and the SSOs are again obtained analytically.

\begin{figure}[!t]
  \centering
  \includegraphics[]{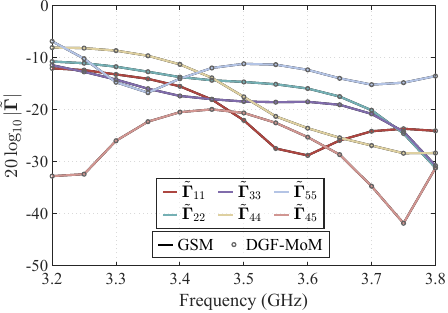}
  \caption{Amplitude of port S-parameters for the multimode horn antenna embedded in a uniaxially anisotropic and homogeneous spherical shell ($\epsilon _r=2,\epsilon _{\bot}=5,\mu _r=1,\mu _{\bot}=3$).}
  \label{fSPara_Horn_singleAnisotropic}
\end{figure}

Commercial solvers such as FEKO, HFSS, and CST do not natively support radially anisotropic material definitions. Nevertheless, an analytical DGF for this medium is available \cite{ref_DGF_anisotropic} and can be incorporated into a MoM solver. Fig.~\ref{fSPara_Horn_singleAnisotropic} compares the S-parameters predicted by our GSM-SSO formulation with those from this specialized DGF-based solver. The excellent agreement confirms the accuracy of the proposed approach for homogeneous anisotropic environments.

\subsection{Radially Piecewise Homogeneous Medium}
\label{SecIV_C}

\begin{figure}[!t]
  \centering
  \includegraphics[]{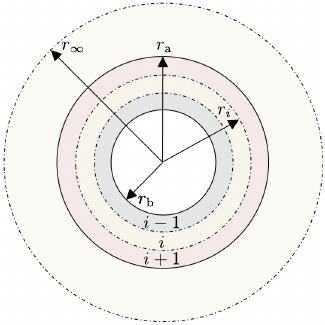}
  \caption{Schematic of a spherically stratified structure consisting of multiple uniformly layered media. The $i$th interface at radius $r = r_i$ separates region $i$ and region $i{+}1$, with region indices increasing outward from the sphere center.}
  \label{fLayered_sph}
\end{figure}

To generalize the previous cases, we now consider a radially piecewise homogeneous medium modeled as a multilayer spherical shell composed of several concentric regions, as illustrated in Fig.~\ref{fLayered_sph}. In this configuration, $\boldsymbol{\epsilon}(r)$ and $\boldsymbol{\mu}(r)$ are piecewise differentiable, and the ODEs in \eqref{eqn16} apply independently within each layer.

Let $g_i(r)$ and $h_i(r)$ be the solutions within the $i$-th layer, following the form of \eqref{eqn18} with the material parameters of that layer. Continuity of tangential fields across $r=r_i$ imposes
\begin{equation}
  \begin{split}
    &g_i\left( r_i \right) =g_{i+1}\left( r_i \right) ,h_i\left( r_i \right) =h_{i+1}\left( r_i \right) \\
    &\frac{g_{i}^{\prime}\left( r_i \right)}{\mu _{\bot i}\left( r_i \right)}=\frac{g_{i+1}^{\prime}\left( r_i \right)}{\mu _{\bot i+1}\left( r_i \right)},\frac{h_{i}^{\prime}\left( r_i \right)}{\epsilon _{\bot i}\left( r_i \right)}=\frac{h_{i+1}^{\prime}\left( r_i \right)}{\epsilon _{\bot i+1}\left( r_i \right)}.
  \end{split}
\end{equation}
These interface conditions enable recursive propagation across layers. For $\mathbf{t}$ and $\mathbf{\Phi}$, the recursion proceeds outward starting from the innermost conditions \eqref{eqn22}, while $\boldsymbol{\rho}$ and $\mathbf{\Psi}$ are obtained by inward recursion beginning with \eqref{eqn29}.

\subsubsection{Two-Layer Isotropic Shell}

\begin{figure}[!t]
  \centering
  \includegraphics[]{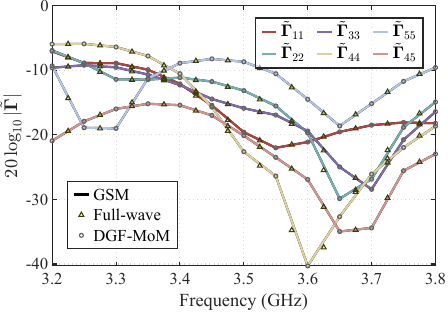}
  \caption{Amplitude of port S-parameters for the multimode horn antenna embedded in a two-layer isotropic and homogeneous spherical shell ($\epsilon_{r1} = \epsilon_{\bot 1} = 4.4 - 0.396\mathrm{j}$, $\epsilon_{r2} = \epsilon_{\bot 2} = 10$). Results from the proposed method, FEKO full-wave simulation, and DGF-MoM.}
  \label{fSPara_Horn_2LayerIsotropic}
\end{figure}

\begin{figure}[!t]
  \centering
  \subfloat[]{\includegraphics[width=0.9in]{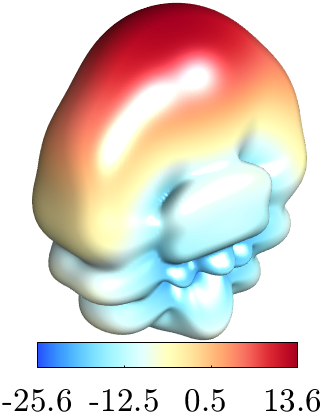}}
  \hfil
  \subfloat[]{\includegraphics[width=0.9in]{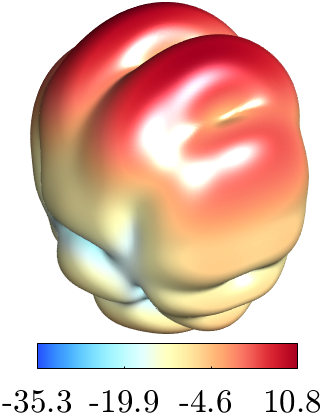}}
  \hfil
  \subfloat[]{\includegraphics[width=0.9in]{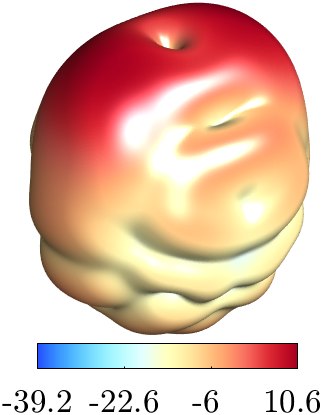}}
  \vfil
  \subfloat[]{\includegraphics[width=0.9in]{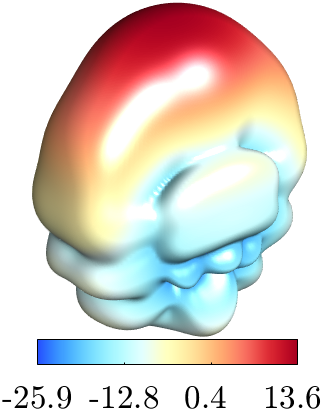}}
  \hfil
  \subfloat[]{\includegraphics[width=0.9in]{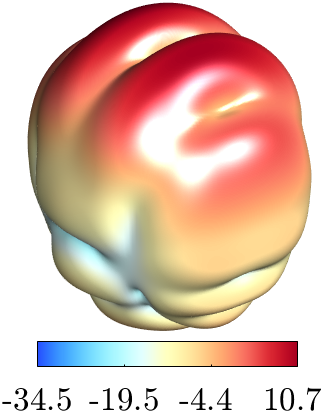}}
  \hfil
  \subfloat[]{\includegraphics[width=0.9in]{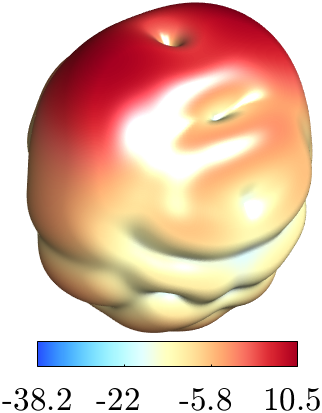}}
  \caption{Far-field gain patterns. The top panel shows results from full-wave simulations using FEKO, while the bottom panel shows results from the proposed method. (a), (d): $\mathrm{TE}_{10}$ excitation; (b), (e): $\mathrm{TE}_{11}$ excitation; (c), (f): $\mathrm{TM}_{11}$ excitation.} 
\label{fGain3D_Horn_2LayerIsotropic} 
\end{figure}

We first consider two isotropic layers of thickness 15~mm, with $\epsilon_{r1}=\epsilon_{\bot 1}=4.4-0.396\j$ and $\epsilon_{r2}=\epsilon_{\bot 2}=10$. Fig.~\ref{fSPara_Horn_2LayerIsotropic} compares the S-parameters obtained using our framework, a DGF–MoM approach, and FEKO. The proposed method matches the DGF-MoM results extremely well; the small discrepancies with FEKO are attributable to meshing of curved surfaces. Far-field patterns in Fig.~\ref{fGain3D_Horn_2LayerIsotropic} also show excellent agreement.

\subsubsection{Two-Layer Anisotropic Shell}

\begin{figure}[!t]
  \centering
  \includegraphics[]{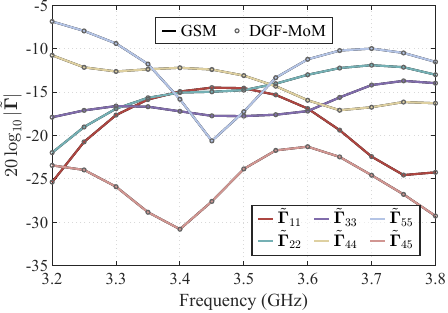}
  \caption{Amplitude of port S-parameters for a multimode horn antenna embedded in a two-layer uniaxially anisotropic and homogeneous spherical shell ($\epsilon _{r1}=2,\epsilon _{\bot 1}=4.4$, $\mu _{r1}=\mu _{\bot 1}=2.2$ and $\epsilon _{r2}=1,\epsilon _{\bot 2}=8$, $\mu _{r2}=2,\mu _{\bot 2}=5$).}
  \label{fSPara_Horn_2LayerAnisotropic}
\end{figure}

We next consider a two-layer uniaxially anisotropic shell with $\epsilon_{r1}=2$, $\epsilon_{\bot 1}=4.4$, $\mu_{r1}=\mu_{\bot 1}=2.2$ and $\epsilon_{r2}=1$, $\epsilon_{\bot 2}=8$, $\mu_{r2}=2$, $\mu_{\bot 2}=5$. Fig.~\ref{fSPara_Horn_2LayerAnisotropic} demonstrates that the proposed method maintains high accuracy for multilayer anisotropic media as well.

\subsection{Radially Piecewise Inhomogeneous Medium}

We finally consider a radially continuous but inhomogeneous medium, \eg
\begin{equation*}
  \begin{split}
    &\epsilon _{\bot 1}\left( r \right) =5\tan \frac{\pi}{5r},\epsilon _{r1}\left( r \right) =1+\exp \left( 2\sin \frac{4}{r} \right) \\
    &\epsilon _{\bot 2}\left( r \right) =2+\ln \left( \frac{2}{r}-5 \right) ,\epsilon _{r2}\left( r \right) =\frac{1}{r}
  \end{split}
\end{equation*}
with $\vec{\mu}(r)=\mat 1_3$ everywhere.

Closed-form solutions of \eqref{eqn16} are not available in this case, so we evaluate $g(r)$ and $h(r)$ numerically using MATLAB's \texttt{ode45}. The equations are non-stiff, and \texttt{ode45} converges efficiently with \texttt{RelTol}~$=10^{-6}$. A reference implementation is provided in \cite{ref_github}.

\begin{figure}[!t]
  \centering
  \includegraphics[]{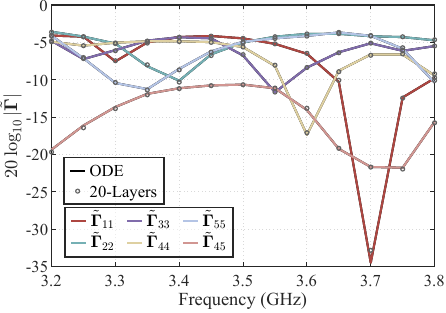}
  \caption{Port S-parameters of the multimode horn antenna embedded in a two-layer radially continuous uniaxially anisotropic spherical shell. The curve labeled ``ODE'' is obtained by numerically solving the radial differential equations \eqref{eqn16}. The circles labeled ``20-Layers'' corresponds to a 20-layer piecewise-homogeneous anisotropic approximation.}
  \label{fSpara_Horn_continue_nonuniform}
\end{figure}

To assess accuracy, we approximate the radial profile using 20 homogeneous anisotropic layers sampled along the profile. Fig.~\ref{fSpara_Horn_continue_nonuniform} compares the ODE-based results with the multilayer approximation. The two solutions agree closely, confirming the correctness of the numerical ODE approach.

\section{Discussion on Different Solving Methods}
\subsection{Applicability and Limitations}

Classical spherical-DGF-based MoM formulations \cite{ref_DGF1,ref_DGF2,ref_DGF3,ref_DGF4,ref_DGF5} exploit analytical layered DGFs so that unknowns are placed only on the antenna surface, guaranteeing rigorous current continuity across metallic junctions even when the structure traverses multiple layers. They impose no constraint on antenna placement and remain the most general framework for multilayer spherical environments. For piecewise-homogeneous media, these schemes deliver high accuracy and efficiency for single-run analyses. Their main limitation is the reliance on closed-form DGFs, which precludes continuously inhomogeneous profiles. In addition, the impedance matrix grows with electrical size and geometric detail and must be rebuilt and refactorized whenever the medium parameters change, making parameter sweeps computationally expensive.

Hybrid T-matrix/MoM schemes \cite{ref_Capek1,ref_Capek2} describe the layered medium via T-matrices (\ie SSOs), thereby avoiding explicit analytical DGFs while still discretizing antenna currents with MoM. In essence, they retain the philosophy of DGF-MoM but implement it in a more numerically flexible manner. A practical constraint, however, is that the antenna must lie entirely within the innermost (or outermost) region, preventing metallic features from crossing layer interfaces. The SSO formulations developed here---including both closed-form operators and ODE-based numerical operators---extend this paradigm to continuously inhomogeneous media and radially uniaxial anisotropy. In this sense, the SSO does not replace the MoM current-domain foundation; rather, it replaces explicit DGF construction with a more general and efficient numerical interface to the medium, extending applicability beyond analytically parameterizable stratifications.

The proposed GSM framework targets implantation-like scenarios where the antenna is fully enclosed by the innermost spherical region. Using spherical waves as global basis functions makes the GSM size depend only on the radius of this region---rather than on antenna geometric complexity---yielding a reduced-order operator independent of radiator detail. Although a MoM formulation is used in earlier sections to extract the antenna GSM, the GSM itself is not tied to MoM; it can, in principle, be constructed by any solver capable of providing port-wave interaction data, as demonstrated later. Once constructed, the GSM combines with the SSOs to account for layered media, including continuously inhomogeneous and anisotropic profiles enabled by the ODE-based SSO. This division of roles---GSM for the radiator and SSO for the surrounding medium---retains full electromagnetic rigor while offering substantial flexibility, enabling geometry-independent compression and reuse of operators, as quantified next.

\subsection{Demonstration on Computational Burden}

\begin{table}[!t]
  \centering
  \caption{Time and memory cost for single evaluation. Executed on an AMD Ryzen 9 9950X 16-Core Processor, 64 GB RAM, 15 threads. The antenna is discretized with 5901 RWGs.}
  \begin{threeparttable}
    \footnotesize 
    \begin{tabular}{@{}lcccc@{}}
      \toprule
      & \multicolumn{2}{c}{\shortstack{$k_{\mathrm f} r_{\mathrm a} = 12$ \\ ($L_{\max}=20$)}}
      & \multicolumn{2}{c}{\shortstack{$k_{\mathrm f} r_{\mathrm a} = 24$ \\ ($L_{\max}=33$)}} \\
      \cmidrule(lr){2-3} \cmidrule(lr){4-5}
      \addlinespace[3pt]
      \textbf{Tasks} & \multicolumn{1}{c}{Time} & \multicolumn{1}{c}{Memory}
                     & \multicolumn{1}{c}{Time} & \multicolumn{1}{c}{Memory} \\
      \midrule
      Full-wave MLFMM      & 77 & 7762 & 231 & 18546 \\
      Full-wave MoM        & 270 & 23860 & $>$8425 & 337646 \\
      \midrule
      \textbf{DGF--MoM}  & & & & \\
      \midrule
      \hspace{1em}\cbullet{} MoM matrix   & 4.2 & 568 & 4.9 & 568 \\
      \hspace{1em}\cbullet{} Solution         & 0.9 & -- & 0.9 & -- \\
      \hspace{1em}\cbullet{} In total         & 5.1 & 568 & 5.8 & 568\\
      \midrule
      \textbf{Proposed Framework} & & & & \\
      \midrule
      \hspace{1em}\cbullet{} Free space MoM     & 3.7 & 308 & 3.7 & 308 \\
      \hspace{1em}\cbullet{} Antenna GSM      & 1.3 & 6 & 1.3 & 6 \\
      \hspace{1em}\cbullet{} SSOs (Ana./Num.) & 0.03/0.1 & $<$0.1/12 & 0.03/0.7 & $<$0.1/12\\
      \hspace{1em}\cbullet{} Solution         & 0.03 & -- & 0.03 & -- \\
      \hspace{1em}\cbullet{} In total         & 5.06/5.13 & 314/326 & 5.06/5.73 & 314/326 \\
      \bottomrule
    \end{tabular}
     \begin{tablenotes}\footnotesize
       \item Ana./Num. denote analytic/numeric evaluations, respectively.
       \item Time in seconds; memory in megabyte (MB).
     \end{tablenotes}
  \end{threeparttable}
  \label{table_case1}
\end{table}

\begin{table}[!t]
  \centering
  \caption{Time and memory cost for single evaluation. Executed on an AMD Ryzen 9 9950X 16-Core Processor, 64 GB RAM, 15 threads. $k_\mrm{f}r_\mrm{a}=12, L_{\max} = 20$.}
  \begin{threeparttable}
    \footnotesize 
    \begin{tabular}{@{}lcccc@{}}
      \toprule
      & \multicolumn{2}{c}{\shortstack{$k_{\mathrm f} R_\mrm{ant.} = 10$ \\ (5901 RWGs)}}
      & \multicolumn{2}{c}{\shortstack{$k_{\mathrm f} R_\mrm{ant.} = 10$ \\ (13273 RWGs)}} \\
      \cmidrule(lr){2-3} \cmidrule(lr){4-5}
      \addlinespace[3pt]
      \textbf{Tasks} & \multicolumn{1}{c}{Time} & \multicolumn{1}{c}{Memory}
                     & \multicolumn{1}{c}{Time} & \multicolumn{1}{c}{Memory} \\
      \midrule
      \textbf{DGF--MoM}   & & & & \\
      \midrule
      \hspace{1em}\cbullet{} MoM matrix   & 4.2 & 568 & 12.7 & 2777 \\
      \hspace{1em}\cbullet{} Solution         & 0.9 & -- & 6.8 & -- \\
      \hspace{1em}\cbullet{} In total         & 5.1 & 568 & 19.5 & 2777\\
      \midrule
      \textbf{Proposed Framework} & & & & \\
      \midrule
      \hspace{1em}\cbullet{} Free space MoM     & 3.7 & 308 & 9.3 & 1433 \\
      \hspace{1em}\cbullet{} Antenna GSM      & 1.3 & 6 & 9.3 & 6 \\
      \hspace{1em}\cbullet{} SSOs (Ana./Num.) & 0.03/0.1 & $<$0.1/12 & 0.03/0.1 & $<$0.1/12\\
      \hspace{1em}\cbullet{} Solution         & 0.03 & -- & 0.03 & -- \\
      \hspace{1em}\cbullet{} In total         & 5.06/5.13 & 314/326 & 18.7/18.7 & 1439/1451 \\
      \bottomrule
    \end{tabular}
  \end{threeparttable}
  \label{table_case2}
\end{table}

Table~\ref{table_case1} compares the proposed GSM framework against full-wave MLFMM, full-wave MoM, and spherical DGF-MoM (hybrid T-matrix/MoM). A homogeneous isotropic shell is used so that all solvers apply without approximation. Two outer electrical radii, $k_\mathrm{f} r_\mathrm{a}=12$ and $k_\mathrm{f} r_\mathrm{a}=24$, are tested. Note that $k_\mathrm{f} r_\mathrm{a}$ measures the electrical size of the outermost spherical layer: while DGF-MoM and the proposed method are dominated by the antenna MoM stage, full-wave solvers scale with both antenna size and spherical-region size. As the spherical domain grows electrically large, full-wave approaches quickly become impractical, whereas DGF-MoM and the proposed GSM formulation remain tractable.

For a fixed medium configuration, the total runtime of the proposed approach is similar to DGF-MoM. The free-space MoM solve and GSM extraction form a one-time cost. Subsequent evaluations require only the compact SSO system: the analytic SSO takes roughly $3\times10^{-2}$~s with sub-megabyte memory, and the ODE-based SSO adds less than 1~s and approximately 12~MB. By contrast, DGF-MoM must rebuild and refactor the full impedance matrix whenever medium parameters change. Once constructed, the GSM occupies only $\sim6$~MB and is reusable for any spherical medium configuration. Thus, although the initial cost is comparable, the proposed method offers orders-of-magnitude savings in both runtime and memory for material sweeps.

Table~\ref{table_case2} examines two antennas sharing the same enclosing radius $R_\mathrm{ant}$ differing in geometric complexity: a simple model with 5,901 RWGs and a more detailed version with 13,273 RWGs. As expected, the MoM matrix size---and therefore the cost of DGF-MoM---increases sharply (from 5.1~s/568~MB to 19.5~s/2777~MB). In the proposed framework, the one-time GSM construction scales with the free-space MoM effort (3.7~s to 9.3~s; GSM extraction 1.3~s to 9.3~s), reflecting geometric complexity. However, the stored GSM itself remains about 6~MB, because its numerical dimension is dictated by the electrical radius $k_\mathrm{f} R_\mathrm{ant}$ rather than by geometric detail. Once available, the same GSM can be paired with the SSOs at negligible cost, demonstrating geometry-independent compression and efficient reuse even for highly complex radiators.

\subsection{Flexibility of the GSM Framework}

This subsection illustrates the flexibility of the proposed GSM framework beyond conventional DGF-MoM formulations. The essential point is that once the free-space GSM of an antenna is available, its behavior inside spherical layered media can be evaluated rapidly using the SSOs developed in this work---regardless of how the GSM was originally obtained.

\subsubsection{Constructing GSM from far-field data}

To demonstrate this idea, we show that the antenna GSM can be assembled directly from radiated far-field patterns (yielding the $\mathbf{R}$ and $\mathbf{T}$ operators) and plane-wave scattering data (yielding the $\mathbf{S}$ matrix). 

Because the vector spherical harmonics form an orthogonal basis, the $ni$-th entry of the transmit matrix $\mathbf{T}$ can be obtained from the radiated far field $\vec F_i^{R}(\hat{\vec r})$ under $i$-th port (eigenmode) single excitation, cf. \eqref{eqn33}:
\begin{equation}
  T_{ni}=\frac{1}{\sqrt{Z_{\mathrm{f}}}}\mathrm{j}^{\tau -l-2}\int_{4\pi}{\vec{A}_n\left( \hat{\vec{r}} \right) \cdot \vec{F}_{i}^{R}\left( \hat{\vec{r}} \right) \mathrm{d}\Omega}
\end{equation}
with the solid-angle measure
\begin{equation*}
  \int_{4\pi}{\mathrm{d}\Omega}=\int_0^{2\pi}{\int_0^{\pi}{\sin ^2\theta \mathrm{d}\theta \mathrm{d}\varphi}}.
\end{equation*}
Reciprocity gives $\mathbf{R}=\mathbf{T}^t$.

Similarly, the $nn^\prime$-th entry of the spherical-wave scattering matrix $\M S$ is obtained from the scattering dyadic kernel $\bar{\mathcal S}(\hat{\vec r},\hat{\vec r}')$:
\begin{equation}
  \begin{split}
      S_{nn^{\prime}}&=\delta _{nn^{\prime}} + 2 \j^{l^{\prime}-l+\tau -\tau ^{\prime}} \times\\
      &\int_{4\pi}{\int_{4\pi}{\vec{A}_n\left( \hat{\vec{r}} \right) \cdot \bar{\mathcal{S}}\left( \hat{\vec{r}},\hat{\vec{r}}^{\prime} \right)}\cdot \vec{A}_{n^{\prime}}\left( \hat{\vec{r}}^{\prime} \right) \mathrm{d}\Omega ^{\prime}}\mathrm{d}\Omega .
  \end{split}
\end{equation}
The dyadic $\bar{\mathcal{S}}\left( \hat{\vec{r}},\hat{\vec{r}}^{\prime} \right) $ relates the vector amplitude $\vec{E}_0$ of an incident plane wave propagating in direction $\hat{\vec r}^\prime$, \ie $\vec E^i(\vec r) = \vec E_0(\hat {\vec r}^\prime)e^{-\j k_\mrm{f} \hat{\vec r}^\prime\cdot \vec r}$, to the corresponding scattered far field in $\hat{\vec r}$ direction:
\begin{equation}
  \vec{F}^S\left( \hat{\vec{r}} \right) =\frac{4\pi}{\mathrm{j}k}\bar{\mathcal{S}}\left( \hat{\vec{r}},\hat{\vec{r}}^{\prime} \right) \boldsymbol{E}_0\left( \hat{\vec{r}}^{\prime} \right)
\end{equation}
see \cite[Eq. (4.29)]{ref_scattering_theory}, \cite{ref_sdyadic_CMA}.

Using these relationships, we have developed a FEKO--MATLAB workflow (released as supplementary material \cite{ref_github}). FEKO serves purely as a free-space electromagnetic engine to generate radiated and scattered far fields; MATLAB then assembles the GSM and couples it with the proposed SSOs to predict antenna performance inside spherical layered media.

Because the required far-field quantities can be supplied by any solver capable of port-driven radiation and plane-wave scattering---not only MoM but also FEM or hybrid FEM--MoM---the methodology extends naturally to a wide range of numerical techniques. Furthermore, since the same fields may also be obtained through antenna measurements, the approach opens a pathway toward experimentally driven GSM-SSO analysis. This enables the inclusion of real-world effects such as fabrication tolerances, connector mismatch, and local perturbations difficult to model numerically. Reliable measurements can also reveal deviations between an actual fabricated radome and its idealized electromagnetic profile, supporting tolerance assessment and statistical evaluation across multiple prototypes.

\begin{figure}[!t]
  \centering
  \includegraphics[width=3in]{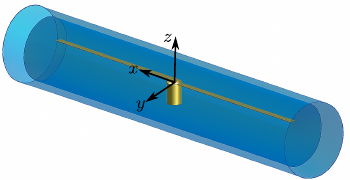}
  \includegraphics[]{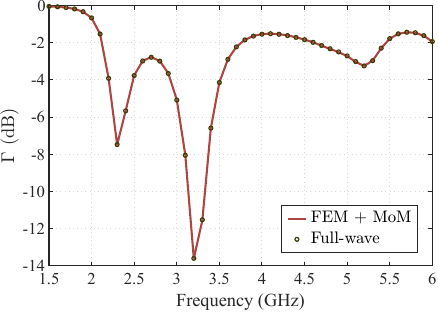}
  \caption{Geometry and performance of the coated dipole antenna. The structure consists of a 50-mm-long, 1-mm-wide dipole enclosed within a cylindrical dielectric layer of relative permittivity 10. The outer dielectric cylinder is 62.5~mm in length and 12~mm in diameter, surrounding an inner cavity of 56.25~mm in length and 10~mm in diameter. The dipole is fed by a 50~$\Omega$ coaxial line (inner/outer radius 0.5/1.15~mm). Also shown is the reflection coefficient when the coated dipole is embedded in a homogeneous spherical shell (inner radius 45~mm, outer radius 60~mm, $\epsilon_r=5-0.5\j$).} 
  \label{fDiople_coat}
\end{figure}

It should be emphasized that the results presented here constitute a numerical demonstration using FEKO as one representative tool. Fig.~\ref{fDiople_coat} shows the case of a dipole enclosed in a cylindrical dielectric coat, analyzed using FEKO's hybrid FEM-MoM solver to generate both radiation and scattering data: metallic parts are modeled by MoM, while dielectrics are treated with FEM. The resulting free-space GSM is then coupled with the SSO framework to predict the port S-parameters inside a dielectric sphere. The predicted results match full-wave simulations extremely well, demonstrating the flexibility and practical applicability of the GSM-SSO methodology for realistic antenna configurations.

\subsubsection{Array-level GSM computation from element data}

Another advantage of the proposed framework is that the GSM naturally supports modular array construction. Array-level GSMs can be synthesized directly from element-level GSMs, enabling rapid evaluation of embedded-element behavior without rerunning full-wave simulations.

This technique was first introduced in \cite{ref_3D_FEM,ref_sph_near_measure} and further developed in \cite{ref_shi_DDM,ref_GSM_overlap,ref_myGSM}. Using the translation formulas of vector spherical waves, the GSM of a single element can be shifted to arbitrary array positions, and mutual coupling among elements can be assembled analytically through the same translation operators. The only required input is the free-space GSM of each isolated element, making the method particularly attractive for array modeling. Because the present framework cleanly separates GSM construction from the SSO stage, a synthesized array-level GSM can immediately be combined with the SSOs to predict array performance inside multilayer spherical environments.

\begin{figure}[!t]
  \centering
  \includegraphics[]{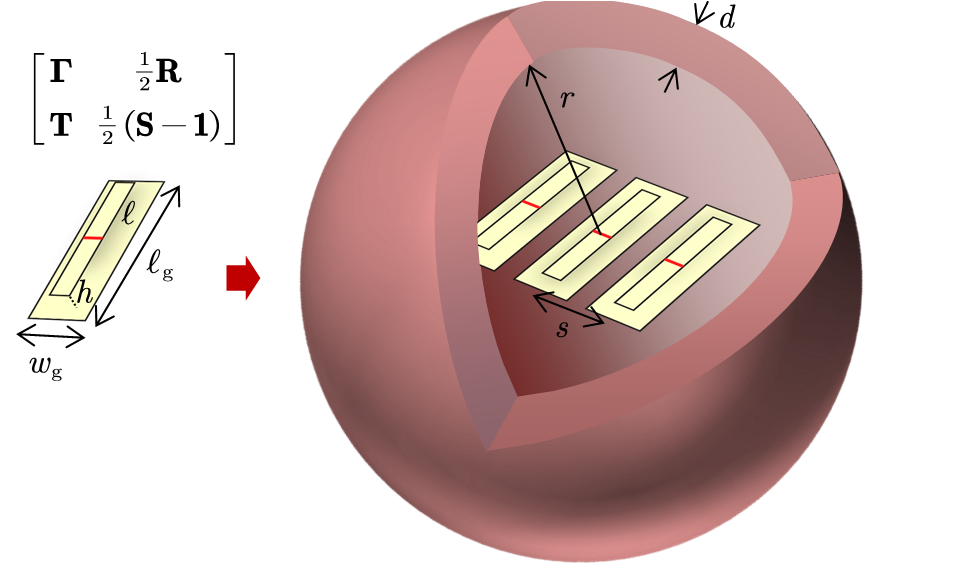}
  \includegraphics[]{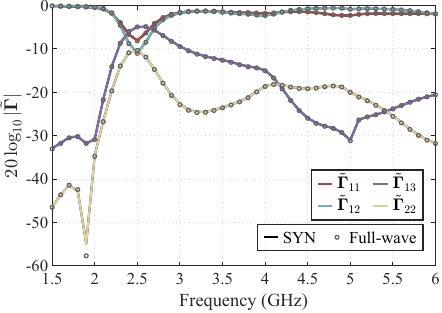}
  \caption{Simulated reflection coefficient of a 3-element dipole array, embedded within a homogeneous spherical shell with relative permittivity $\epsilon_r=5-0.5\j$. The dipole dimensions are $\ell=52$~mm, $w=0.1\ell$, $\ell_g=1.25\ell$, $w_g =0.25\ell$, and $h=0.31\ell$. The spherical-shell radius and thickness are $r=50$~mm and $d=10$~mm. Element spacing is $s=20$~mm.}
  \label{fSpara_Diople_MoM_FEM}
\end{figure}

Figure~\ref{fSpara_Diople_MoM_FEM} demonstrates a simple example. The free-space GSM of a dipole-strip element with a reflector is first extracted. Following the procedure in \cite[Section~II]{ref_shi_DDM}, the free-space GSM of a three-element array is synthesized from the individual element-level GSM. This synthesized GSM is then combined with the proposed GSM-SSO framework to evaluate the port characteristics of the array when placed inside a dielectric spherical shell (``SYN'' curve in Fig.~\ref{fSpara_Diople_MoM_FEM}). The excellent agreement with full-wave simulations confirms the practicality of the approach and highlights the potential of modular GSM assembly for embedded array analysis.

\section{Conclusion}

This work has presented a unified framework for modeling antennas embedded in spherically layered media. The approach is built on the free-space GSM of the antenna and a set of semi-analytical SSOs describing the surrounding medium. Because the GSM depends only on the radius of the enclosing sphere rather than on antenna geometry, the method achieves inherent model-order reduction. As a result, once the free-space GSM is obtained, antenna behavior under arbitrary modifications of the spherical medium can be evaluated without re-simulating the radiator, yielding orders-of-magnitude improvements in efficiency over classical DGF-MoM schemes.

The framework is solver-agnostic. We demonstrated that the GSM can be constructed directly from radiated and scattered far-field data, enabling its acquisition from various numerical solvers—and, in principle, from measurements. This provides a natural pathway for incorporating real-world effects such as fabrication tolerances and material deviations.

The GSM formulation is also modular: array-level GSMs can be synthesized analytically from element-level GSMs via spherical-wave translation, offering rapid evaluation of embedded-element and array performance inside layered spherical media without rerunning array-scale full-wave models. Examples confirm that the synthesized array GSM, combined with the SSO operators, accurately reproduces full-wave results.

Finally, the SSO derivations and accompanying ODE-based numerical evaluations extend the framework to radially inhomogeneous and radially uniaxial anisotropic structures. Although fully anisotropic, continuously varying media remain an open challenge, this affects only SSO construction and does not compromise the generality of the GSM-SSO paradigm.
\newpage

\begin{appendices}
\section{Vector Spherical Waves and Special Functions}
\label{APP_A}
Vector spherical harmonics serve as essential building blocks for expanding vector fields in spherical coordinates. They are defined as follows \cite[Appendix C.4]{ref_scattering_theory}:
\begin{equation}
  \begin{split}
    &\vec{A}_{1\sigma ml}\left( \hat{\vec{r}} \right) =\frac{1}{\sqrt{l\left( l+1 \right)}}\nabla Y_{\sigma ml}\left( \hat{\vec{r}} \right) \times \hat{\vec{r}}\\
    &\vec{A}_{2\sigma ml}\left( \hat{\vec{r}} \right) =\frac{1}{\sqrt{l\left( l+1 \right)}}r\nabla Y_{\sigma ml}\left( \hat{\vec{r}} \right) \\
    &\vec{A}_{3\sigma ml}\left( \hat{\vec{r}} \right) =\hat{\vec{r}}Y_{\sigma ml}\left( \hat{\vec{r}} \right) .
  \end{split}
\end{equation}
Here, $Y_{\sigma ml}\left( \hat{\vec{r}} \right) $ denotes the scalar spherical harmonic, given by
\begin{equation}
  Y_{\sigma lm}\left( \hat{\vec{r}} \right) =\sqrt{\frac{2-\delta _{m0}}{2\pi}}\tilde{P}_{l}^{m}\left( \cos \theta \right) \begin{cases}
    \cos \left( m\varphi \right)\\
    \sin \left( m\varphi \right)\\
  \end{cases}
\end{equation}
where $\tilde{P}_l^m(x)$ is the normalized associated Legendre function and $\delta_{ij}$ is the Kronecker delta.

Using these harmonics, the vector spherical wave functions can be expressed compactly as:
\begin{equation*}
  \begin{split}
    &\vec{u}_{1\sigma ml}^{\left( p \right)}\left( k\vec{r} \right) =R _{1l}^{\left( p \right)}\left( kr \right) \vec{A}_{1\sigma ml}\left( \hat{\vec{r}} \right) \\
    &\vec{u}_{2\sigma ml}^{\left( p \right)}\left( k\vec{r} \right) =R _{2l}^{\left( p \right)}\left( kr \right) \vec{A}_{2\sigma ml}\left( \hat{\vec{r}} \right) +R _{3l}^{\left( p \right)}\left( kr \right) \vec{A}_{3\sigma ml}\left( \hat{\vec{r}} \right) .
  \end{split}
\end{equation*}

The scalar radial functions $R_{\tau l}^{(p)}(kr)$ are constructed from Riccati-type special functions:
\begin{equation}
  \begin{split}
    &R _{1l}^{\left( p \right)}\left( kr \right) =\frac{z_{l}^{\left( p \right)}\left( kr \right)}{kr}     \\
    &R _{2l}^{\left( p \right)}\left( kr \right) =\frac{\left[ z_{l}^{\left( p \right)}\left( kr \right) \right] '}{kr}    \\
    &R _{3l}^{\left( p \right)}\left( kr \right) =\sqrt{l\left( l+1 \right)}\frac{z_{l}^{\left( p \right)}\left( kr \right)}{(kr)^2}
  \end{split}
\end{equation}
where $z_l^{(1)}(x) \equiv \zeta_l(x)$ and $z_l^{(4)}(x) \equiv \xi_l(x)$ are the Riccati–Bessel and Riccati–Hankel functions, respectively. Their explicit definitions are:
\begin{equation}
  \begin{split}
    \zeta_l(x)=\sqrt{\frac{\pi x}{2}}J_{l+\frac{1}{2}}(x)
    \xi_l(x)=\sqrt{\frac{\pi x}{2}}H^{(2)}_{l+\frac{1}{2}}(x)
  \end{split}
\end{equation}
with $J_n(x)$ and $H_n^{(2)}(x)$ denoting the cylindrical Bessel and second-type Hankel functions, respectively.

\end{appendices}

\begin{IEEEbiography}[{\includegraphics[width=1in,height=1.25in,clip,keepaspectratio]{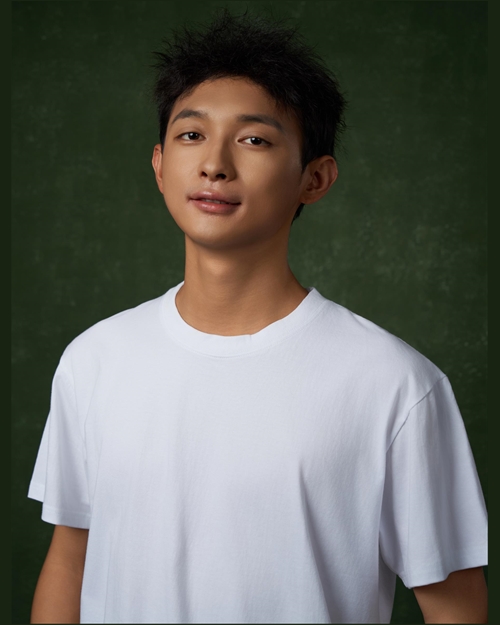}}]{Chenbo Shi}
  Chenbo Shi was born in 2000 in China. He received his Bachelor's degree from the University of Electronic Science and Technology of China (UESTC) in 2022. He is currently pursuing his Ph.D. at the same institution. His research interests include electromagnetic theory, characteristic mode theory, and computational electromagnetics. 
  
  Chenbo has been actively involved in several research projects and has contributed to publications in these areas. His work aims to advance the understanding and application of electromagnetic phenomena in various technological fields.
  \end{IEEEbiography}

  \begin{IEEEbiography}[{\includegraphics[width=1in,height=1.25in,clip,keepaspectratio]{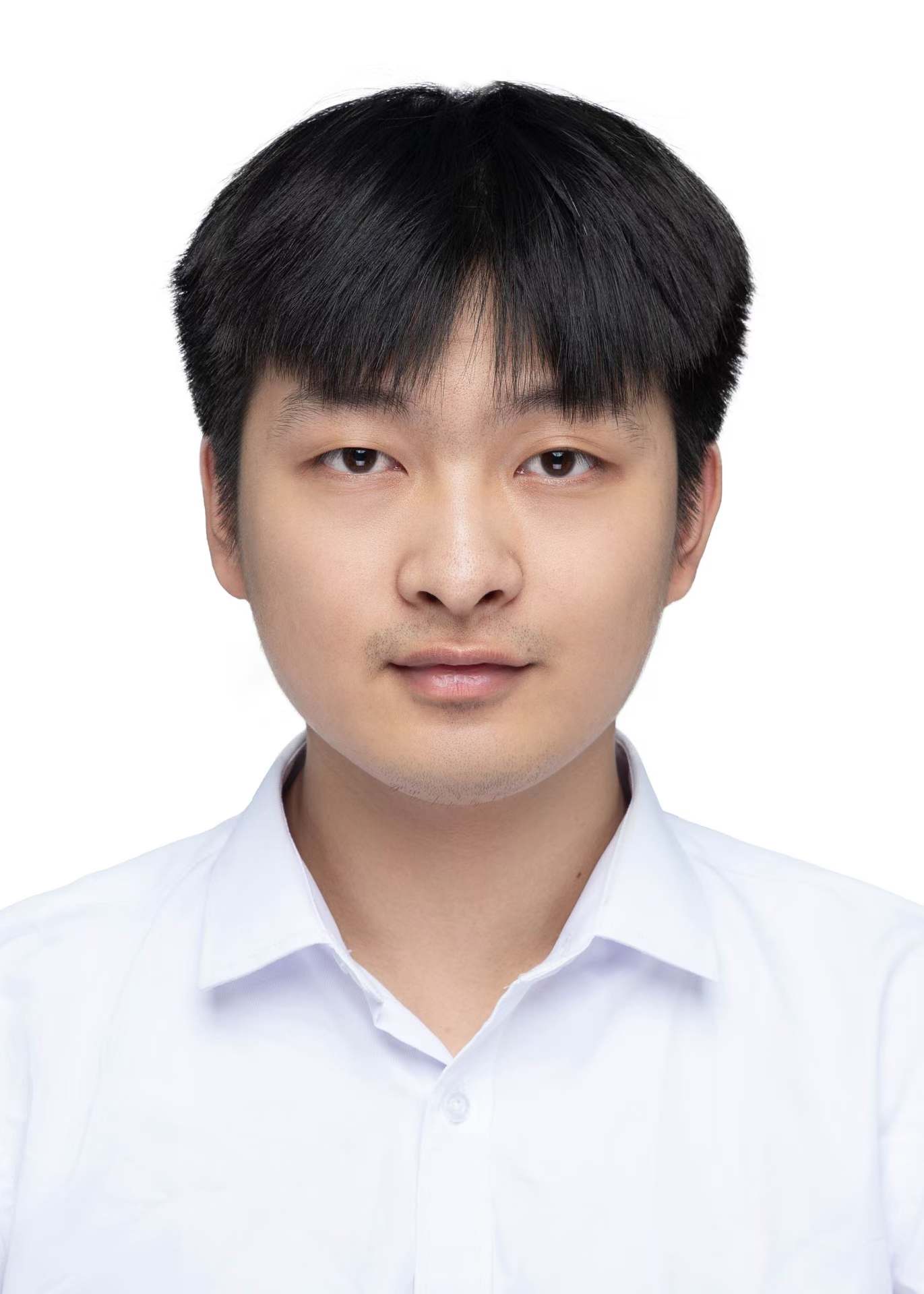}}]{Xin Gu}
  received the B.E. degree from Chongqing University of Posts and Telecommunications (CQUPT), ChongQing, China, in 2022. He is currently pursuing the Ph.~D degree with the School of Electronic Science and Engineering, University of Electronic Science and Technology of China (UESTC), Chengdu, China.
    
  His research interests include electromagnetic theory and electromagnetic measurement techniques
  \end{IEEEbiography}

  \begin{IEEEbiography}[{\includegraphics[width=1in,height=1.25in,clip,keepaspectratio]{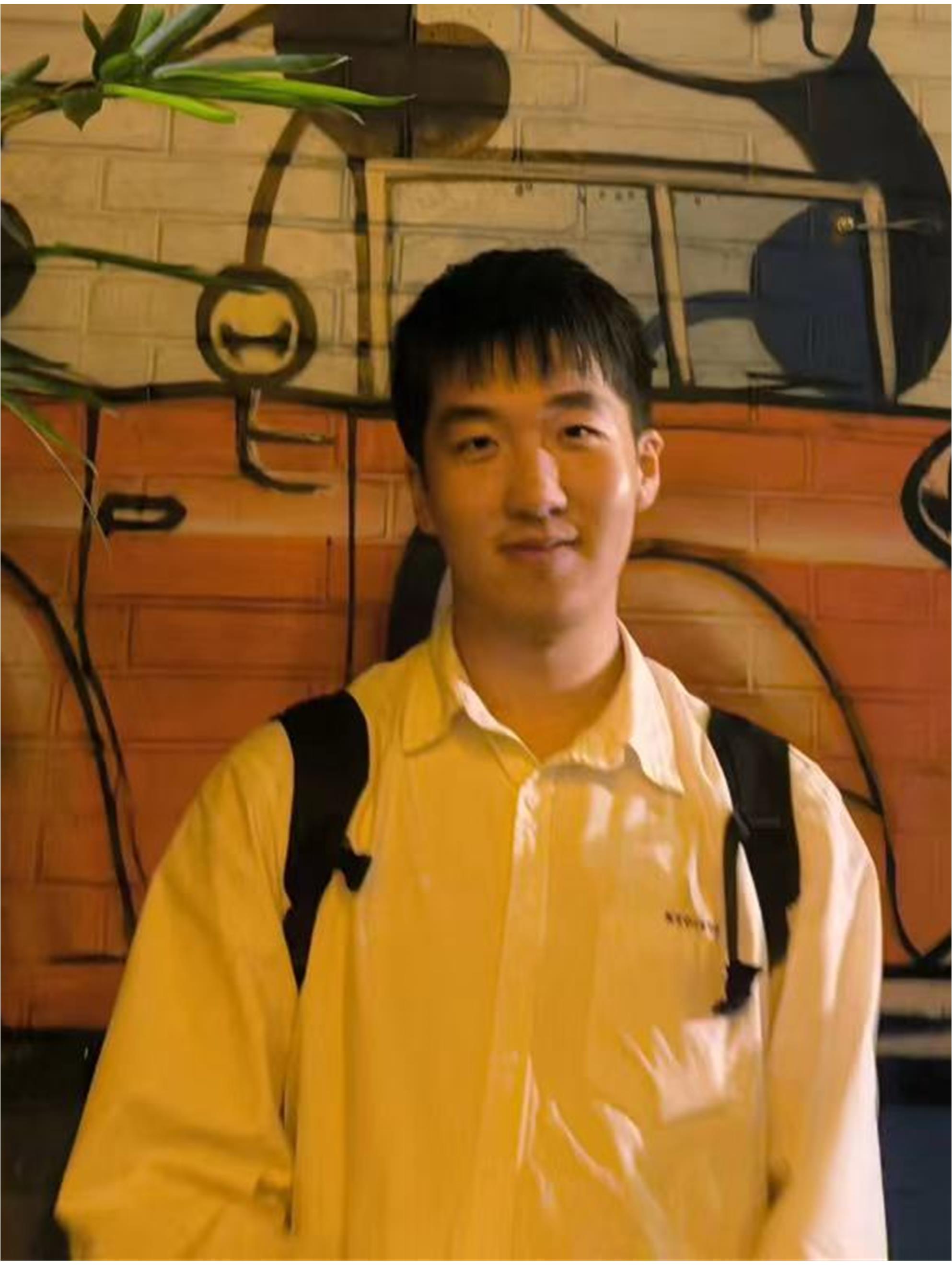}}]{Shichen Liang}
  received the B.E. degree from Beijing University of Chemical Technology (BUCT), Beijing, China, in 2022. He is currently pursuing the Ph.~D degree with the School of Electronic Science and Engineering, University of Electronic Science and Technology of China (UESTC), Chengdu, China. 
  
  His research interests include electromagnetic theory and electromagnetic measurement techniques.
  \end{IEEEbiography}

 \begin{IEEEbiography}[{\includegraphics[width=1in,height=1.25in,clip,keepaspectratio]{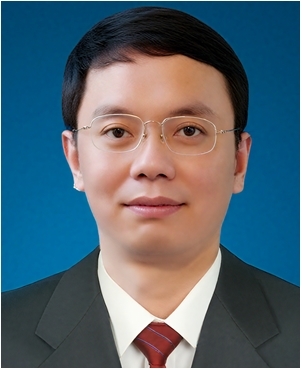}}]{Jin Pan}
    received the B.S. degree in electronics and communication engineering from the Radio Engineering Department, Sichuan University, Chengdu, China, in 1983, and the M.S. and Ph.D. degrees in electromagnetic field and microwave technique from the University of Electronic Science and Technology of China (UESTC), Chengdu, in 1983 and 1986, respectively. 
    
    From 2000 to 2001, he was a Visiting Scholar in electronics and communication engineering with the Radio Engineering Department, City University of Hong Kong. He is currently a Full Professor with the School of Electronic Engineering, UESTC. 
    
    His current research interests include electromagnetic theories and computations, antenna theories, and techniques, field and wave in inhomogeneous media, and microwave remote sensing theories and its applications. 
    \end{IEEEbiography}

\end{document}